\documentclass[a4paper]{amsart}

\newcommand*{\mykeywords}{Exponential rings, exponential polynomials, exponential ideals}

\title{E-ideals in exponential polynomial rings}

\usepackage{lmodern}
\usepackage{microtype}

\author[P. D'Aquino]{Paola D'aquino}
\address{Universit\`a della Campania "L. Vanvitelli"}
\email{paola.daquino@unicampania.it}

\author[A. Fornasiero]{Antongiulio Fornasiero}
\address{Universit\`a di Firenze}  
\email{antongiulio.fornasiero@gmail.com} 
\urladdr{https://sites.google.com/site/antongiuliofornasiero/}

\author[G. Terzo]{Giuseppina Terzo}
\address{Universit\`a degli Studi di Napoli "Federico II"}
\email{giuseppina.terzo@unina.it}

\thanks{Members of the ``National Group for Algebraic and Geometric Structures,
  and their Applications'' (GNSAGA - INdAM). 
This research is part of the project  PRIN 2022 ``Models, sets and classification''.}

\usepackage{ifpdf}
\ifpdf
\usepackage[pdftex,
pdfborder={0 0 0}, plainpages=false,%
pdfpagelabels, pdfdisplaydoctitle=true,%
pdfstartview={XYZ null null null},%
colorlinks=false,
pdfusetitle
]{hyperref}
\hypersetup{%
 pdfauthor={P. D'Aquino, A. Fornasiero, G. Terzo},
  pdfkeywords={\mykeywords},
  pdfsubject={Exponential ideals in the ring of exponential polynomials:
    maximal and prime}
}
\else
\usepackage[plainpages=false, linktocpage=true, pdfpagelabels, colorlinks=false
]{hyperref}
\fi

\usepackage{amsmath, amsthm, amssymb}
\usepackage{mathtools}
\usepackage{xspace}

\usepackage{comment}

\usepackage[neveradjust]{paralist}
\AtBeginDocument{\setdefaultleftmargin{0pt}{}{}{}{}{}}
\setdefaultenum{\upshape(1)}{\upshape i)}{}{}

\usepackage[T1]{fontenc}

\usepackage{url}

\usepackage[msc-links]{amsrefs}

\makeatletter
\@namedef{subjclassname@2020}{%
 \textup{2020} Mathematics Subject Classification}
\makeatother

\newcommand*{\intro}[1]{\textbf{#1}}
\newcommand*{\Pa}[1]{\bigl( #1 \bigr)}
\newcommand*{\set}[1]{\{#1\}}

\newcommand*{\card}[1]{\lvert#1\rvert}
\newcommand{\N}{\mathbb{N}}
\newcommand{\Z}{\mathbb{Z}}
\newcommand{\R}{\mathbb{R}}
\newcommand{\C}{\mathbb{C}}
\newcommand{\Q}{\mathbb{Q}}

\newcommand*{\tuple}[1]{\langle #1 \rangle}

\newcommand{\tv}{\bar t}
\newcommand{\av}{\bar a}

\newcommand{\cv}{\bar c}

\newcommand{\x}{\bar x}
\newcommand{\y}{\bar y}

\newcommand{\Ffam}{\mathcal F}
\newcommand{\Tbar}{\overline T}
\DeclareMathOperator{\Erad}{E-rad}
\newcommand*{\model}[1]{\mathcal M(#1)}
\newcommand*{\radicale}[1]{\operatorname{{#1}-rad}}
\DeclareMathOperator{\Tzprimerad}{T^\prime_0-rad}
\DeclareMathOperator{\Tradop}{T-rad}
\DeclareMathOperator{\Tzradop}{T_p-rad}
\DeclareMathOperator{\Tunoradop}{T_1-rad}
\newcommand*{\Trad}[1]{\Tradop(#1)}
\newcommand*{\Tzrad}[1]{\Tzradop(#1)}
\newcommand*{\Tunorad}[1]{\Tunoradop(#1)}

\def\Ind#1#2{#1\setbox0=\hbox{$#1x$}\kern\wd0\hbox to
  0pt{\hss$#1\mid$\hss}\lower.9\ht0\hbox to 0pt{\hss$#1\smile$\hss}\kern\wd0}

\newcommand{\wloG}{w.l.o.g\mbox{.}\xspace}

\newcommand{\ie}{i.e\mbox{.}\xspace}

\newtheorem{lemma}{Lemma}[section]
\newtheorem{thm}[lemma]{Theorem}
\newtheorem{corollary}[lemma]{Corollary}

\newtheorem{proposition}[lemma]{Proposition}
\newtheorem{open problem}[lemma]{Open problem}

\newtheorem*{fact*}{Fact}

\theoremstyle{remark}

\newtheorem{claim}{Claim}
\newtheorem*{claim*}{Claim}

\newtheorem{remark}[lemma]{Remark}

\theoremstyle{definition}
\newtheorem{definition}[lemma]{Definition}
\newtheorem{notation}[lemma]{Notation}
\newtheorem{final remark}[lemma]{Final remark}
\newtheorem{example}[lemma]{Example}
\newtheorem{examples}[lemma]{Examples}

\begin{document}

\begin{abstract}
We investigate exponential ideals within the context of exponential polynomial
rings over exponential fields. We establish two distinct notions of maximality
for exponential ideals and explore their relationship to primeness. 

These three concepts--prime, maximal, and E-maximal--are 
shown to be independent, in contrast to the classical scenario.
Furthermore, we  demonstrate that, over an algebraically
closed field K, the correspondence between points of $K^n$ and maximal exponential ideals of the ring of exponential polynomials breaks down. Finally, we introduce and characterize exponential radical ideals.
\end{abstract}

\maketitle


\makeatletter
\renewcommand\@makefnmark%
   {\normalfont(\@textsuperscript{\normalfont\@thefnmark})}
\renewcommand\@makefntext[1]%
   {\noindent\makebox[1.8em][r]{\@makefnmark\ }#1}
\makeatother

\section{Introduction}

An exponential ring (or E-ring for short)
is  a commutative unital ring $R$ together with a monoid homomorphism $E: (R, 0,
+) \to (R^*, 1, \cdot)$.
Notable examples include the real and complex fields equipped with their standard exponential functions.

Given an E-ring $R$ we consider  exponential polynomials over $R$ in the variables $\overline{x}$, i.e.\ formal expressions constructed from constants in $R$ and $\overline{x}$ using $+$, $\cdot$ and E, (see \S2 for definitions and details of the construction). In a natural way an E-ring structure can be defined over the exponential polynomials over $R$, and we denote this ring by $R[\overline x]^E$. 
The ring of exponential polynomials with only one iteration of exponentiation
over the complex field has been studied intensively from an analytic point of
view since the early 1920s by Polya, Ritt, Schwengeler, and others in order to
understand the distribution of their roots. It is worth recalling the work of
Henson and Rubel \cite{HR}, and Katzberg \cite{Katzberg} where they
characterize  the exponential polynomials over $\C$ with no roots, and those
with finitely many roots, respectively.

The model theoretic analysis of exponential fields starts with the problem posed
by Tarski in the 30's on the decidability of the field of reals with the usual
exponential function. Many authors have contributed towards a solution of this
problem in over 60 years. In the late 90's  Macintyre and Wilkie gave a positive
answer to Tarski's problem assuming Schanuel's conjecture in transcendental
number theory,   see \cite{macWilkie}. In a more general setting,  exponential
rings have been analyzed in \cites{macintyre, van}, just to mention a few articles.

 In recent years, a new impulse to the study of exponential structures, and in
 particular of exponential polynomials,  originated from work of Zilber on a
 possible axiomatization of the complex exponential field. Some axioms have been
 formulated
in terms of zero sets of systems of exponential polynomials. Using
 Zilber's axioms and purely algebraic and geometrical  arguments, properties
 known for exponential polynomials over $\mathbb C$ have been extended to
 exponential polynomials and systems of them over the E-fields identified by
 Zilber, see for example \cites{null, shapiro}. Vice versa, using methods
 and results from Diophantine geometry many instances of one of Zilber's axioms
 have been shown to be true for
 exponential polynomials over $\mathbb C$, see \cites{marker, mantova, bayskirby, DFT1, DFT2}. 
Further contributions in this direction can be found in the very recent publications  \cite{Arsla} and \cite{Gall}.

Given an exponential ring $R$,
an exponential ideal (or E-ideal for short) is an ideal $I$ of $R$ such that, for every $a \in I$, $E(a) - 1 \in I$; 
  equivalently, an E-ideal is the kernel of a homomorphism of exponential rings with domain $R$ (see \S\ref{sec:E-ideal} and \cite{macintyre}). A source of inspiration for this work is differential algebra, where people study (among other things) differential ideals of the ring of differential polynomials $\C\{X\}.$
In that context, it is known that maximal differential ideals are prime. Moreover, one of the most important result is Ritt-Raudenbush Theorem, asserting that there is no infinite ascending chain of radical differential ideals (and therefore $\C^n$ with the Kolchin topology is a Noetherian space) \cite{Kolchin}.

We show that for exponential algebra the situation is quite different. 
As Macintyre observed, some of the fundamental  results true in polynomial rings
do not hold in the E-rings  of exponential polynomials, e.g., Hilbert Basis
Theorem. Indeed, the E-ring of exponential polynomials over $\C$ (and also over
$\R$) is not Noetherian since the E-ideal 
$I = (E(\frac{x}{2^n}) - 1)_{n \in \N}^E$ 
is not finitely generated, see also \cites{MAC_free, HRS, terzo}.  
As stated in  \cite{macintyre} $(\C, E)$  is not Noetherian with respect to  the exponential Zariski topology where the basic closed sets are the zero sets of exponential polynomials over $\C$. On the contrary, $(\R, E)$ is Noetherian as topological space in the exponential Zariski topology (see~\cite{tougeron}).

To our knowledge, a comprehensive and systematic analysis of E-ideals remains lacking.  
In \cites{MAC_free, terzo}  E-ideals associated to suitable  E-ring homomorphisms  have been characterized in order to  
identify the algebraic relations between $\pi$ and $i$ in $(\mathbb C, E)$, and for determining the E-subring of $(\mathbb R, E)$ generated by $\pi$.  In both cases the characterizations are obtained modulo Schanuel's conjecture. Recently, Point and Regnault in \cite{point} made some progress towards a version of a Nullstellensatz for exponential polynomials over an E-ring. In this paper we improve and extend some of their results.

We consider the following three notions for E-ideals:
\begin{enumerate}
\item prime E-ideals;
\item E-ideals which are maximal as ideals, which we  call "strongly maximal";
\item E-ideals which are maximal among E-ideals, which we  call "E-maximal".
\end{enumerate}

We want to study these three categories of E-ideals inside $R[\overline x]^E$, where $(R,E)$ is an E-ring. In order to simplify the exposition in this introduction we  consider the case when $R = \C$, and we show that some of the fundamental properties of classical ideals of $\C[x]$ do not generalize to E-ideals of $\C[x]^E$.

Obviously, strongly maximal ideals are both prime and E-maximal.
We show that the converse does not hold.

\begin{thm}
There exists an E-maximal E-ideal of $\C[x]^E$ which is not prime
(see Lemma~\ref{maxnoprime}).

There exists a nonzero prime E-ideal which is not E-maximal
(see Example~\ref{ex:prime-not-max}).

There exists a prime and E-maximal E-ideal which is not strongly maximal
(see Lemma~\ref{lem:emax-not-max}).

There exists a prime E-ideal which is not the intersection of E-maximal
E-ideals (see Lemma~\ref{lem:notHilbert}).
\end{thm}

On the other hand, we have
\begin{thm}
If $f \in \C[x]^E$ is irreducible, then the E-ideal $(f)^E$ generated by $f$ is
prime (see Proposition~\ref{E-prime}).

Moreover, for any given $a\in \C$ the set 
\[
I_{a} := \{f \in \C[x]^E: f(a) = 0\}
\]
is a strongly maximal E-ideal (see Lemma~\ref{maximal}).
\end{thm}

\medskip


One of the difficulties in the study of the ring of exponential polynomials is
that, unlike the case of classical polynomials or of differential polynomials, there is no
``universal domain'' where to look for zeros of exponential polynomials.
Moreover, as suggested by Angus Macintyre, we prove that the weak Nullstellensatz fails in our context.

\begin{thm}
There exists a strongly maximal E-ideal of $\C[x]^E$ with no zeroes in~$\C$
(see Proposition~\ref{thm:max-nontrivial} and Theorem~\ref{thm:max-zero}).
The exponential ring $\C[x]^E$  does not satisfy the weak Nullstellensatz (see Proposition~\ref{Null}).
\end{thm}

\smallskip

In \S\ref{sec:Ekernel} we restrict ourselves to the study of E-ideals which
 do not
``increase'' the kernel of the exponential function.
We show that a \emph{prime} E-ideal which does not increase the kernel may still
have no zeroes in $\C$ (Corollary~\ref{cor:E-kernel}): we leave open  what happens for maximal E-ideals.

\medskip

In order to study E-ideals of $\C[x]^E$, we introduce E-ideals of \emph{partial}
E-rings (see Definition \ref{partial}), and we describe a filtration of $\C[x]^E$ by partial
E-rings (see  \S\ref{partialrings} and  \S\ref{sec:E-ideal}), expanding the work
of \cites{van, macintyre, Manders, point}.  
One reason for this approach is that E-ideals of $\C[x]^E$ are difficult
to build and study. For instance, it is not clear how to prove directly that a
certain E-ideal (given by its generators) is proper and prime.  We consider
instead a much smaller ring S of $\C[x]^E$ (e.g., $S = \C[x, e^{\Q x}]$, as in the proof of
Theorem \ref{maxnoprime}), which is only a partial E-ring. We can build rather easily a
(proper) prime E-ideal $I$ of S, and we use two general results (Lemma \ref{lem:R-S} and
Corollary \ref{lem:ideal-1-extension}) to show that, under suitable assumptions
on S, $I$ generates a prime E-ideal of $\C[x]^E$.  
Similarly, in Proposition 5.4, in order to show that an irreducible element of $\C[x]^E$ generates a prime E-ideal, we work first in a smaller partial E-subring of $\C[x]^E$. This strategy works also in some of our other examples.

\medskip

We also consider a fourth notion: what should be the ``right'' definition of
E-radical E-ideal?
We show that simply saying that ``$J$ is an E-ideal  which is radical as an ideal''
is not a satisfactory definition (see Remark \ref{idealeradicale}: the E-ideal of $\C[x,y]^{E}$
generated by $xy$  is a radical ideal which is not an intersection of prime E-ideals).
We propose the following:
\begin{definition} 
An E-ideal is E-radical if it is an intersection of prime E-ideals.

The E-radical of an E-ideal $J$ is the intersection of the prime E-ideals
containing~$J$.
\end{definition}

We give equivalent characterizations of E-radical ideals and alternative ways to
build the E-radical of an E-ideal (Theorems~\ref{thm:Erad} and
\ref{thm:Erad-up}), using among others tools from universal algebra.

\section{Partial E-domains and E-polynomial rings}\label{partialrings}

In this section we recall the relevant notions relative to exponential  rings and exponential polynomial rings. All rings we  consider have characteristic $0$.

\begin{definition}\label{partial}
A partial E-ring (or partial exponential ring) is a triple $D=(D; V,E)$ where
\begin{enumerate}
\item $D$ is a   $\Q$-algebra; 
\item $V $ is a  $\Q$-vector subspace of
  $D$ containing~$\Q$;
\footnote{In the most interesting cases, $V$ is a  $\Q$-subalgebra of~$D$.}
\item $E: (V, +) \to (D^{*}, \cdot)$ is a group homomorphism.
\end{enumerate}
When $D$ is an integral domain, we refer to  $(D; V,E)$ as a partial  E-domain. If $V=D$ then the exponential function is total on $D$, and we  call it simply an E-ring (or exponential ring). 
\end{definition}

An extension of a partial E-ring $(D; V, E)$ is given by a partial E-ring $(D'; V', E')$ and an injective map $\phi : D \to D'$ such that: $\phi$ is a homomorphism of $\Q$-algebras, and, for every $v \in V$, $\phi(v) \in V'$ and $E'(\phi(v))  =  \phi(E(v)).$\\



 We fix  a partial E-domain $D = (D; V, E)$, and we describe how to extend formally the exponential function in $D$.\\

\begin{definition}
A partial free 1-extension of $D$ is a partial E-domain constructed in the
following way.

Let $A$ be a $\Q$-vector subspace of $D$ such that $A \cap V = \{0\}$, and 
consider the group ring $D[t^{A}]$ as a $D$-algebra via $t^a\cdot t^b=t^{a+b}$. Notice that $A$ is torsion free, and so $D[t^{A}]$ is a domain.
Now we extend the partial exponential map to $A$ 
with values in $D[t^{A}]$. 
Let $V' = V \oplus A$ be the $\Q$-vector subspace of $D$. 
Any  element of $V'$ can be written uniquely as $v + a$, for some $v \in V$
and $a \in A$. We  define the function $E': (V',+) \to (D[t^{A}]^*,\cdot)$ as follows: 
for any given $v \in V$ and $a \in A$, we  set
\[
E'(v + a) := E(v) t^{a}.
\]
Then, $D' = (D[t^{A}],V',E')$ is a partial free 1-extension of~$D$.

\end{definition}
Notice that $D'$ is indeed a partial E-domain.

\begin{definition}
\label{def:1-ext}
If $A$ is a complement of $V$ inside $D$ as $\Q$-vector spaces, we say
 that $D[t^{A}]$ is the free 1-extension of $D$, (and we denote it by $\bf D'$).
\end{definition}

\begin{lemma}
\label{unico}
1) The free 1-extension  ${\bf D}'=D[t^{A}]$ of $D$ is unique: that is, it does not depend on the choice of~$A$.
In other words, if $B $ is another complement of $V$ in $D$, then $D[t^B]$ and $D[t^{A}]$ are isomorphic (over $D$) as partial E-domains.

2) The free 1-extension ${\bf D}'$ is the universal free extension of $D$ having $D$ as domain of the
exponential map.
In other words, if $F$ is any partial E-ring extending $D$ and such that $D$
is contained in the domain of the exponential function of $F$, then there is a
unique homomorphism of E-domains over $D$ from ${\bf D}'$ into $F$.
\end{lemma}
\begin{proof}
1) Let $\alpha: D \to A$ and $\upsilon: D \to V$ be the maps corresponding to the decomposition
$D = V \oplus A$: i.e., for any $d\in D$, $\alpha (d)\in A$ and $\upsilon (d)\in V$ are the unique
elements of $A$ and $V$ respectively, such that $d= \upsilon (d)+\alpha (d)$. 
Define the map $\psi: t^{B} \to D[t^{A}]$ by
\[
\psi(t^{b}) := E(\upsilon(b)) t^{\alpha(b)}
\]
and extend $\psi$ by $D$-linearity to $D[t^{B}]$.
Then, $\psi$ is an isomorphism of partial E-domains fixing $D$ pointwise.
In fact, it is clear that $\psi$ is an isomorphism of additive groups and that it
fixes $D$ pointwise.
Given $x, x' \in D[t^{B}]$, write
\[
x = \sum_{b \in B} d_{b} t^{b} \qquad  x' = \sum_{b' \in B} e_{b'} t^{b'} 
\]
for some (unique) $d_{b}, e_{b'} \in D$.
Then,
\begin{multline*}
\psi(x x') 
= \psi(\sum_{b,b' \in B} d_{b} e_{b'} t^{b + b'} )
= \sum_{b,b' \in B} d_{b} e_{b'} E(\nu(b + b')) t^{\alpha(b + b')}
= \\ 
= \sum_{b,b' \in B} d_{b} e_{b'} E(\nu(b) + \nu(b')) t^{\alpha(b)+ \alpha(b')}
= \\
= \sum_{b,b' \in B} d_{b} e_{b'} E(\nu(b)) E(\nu(b')) t^{\alpha(b)} t^{\alpha(b')}
= \psi(x) \psi(x'),
\end{multline*}
where we use the fact that $\nu(b + b') = \nu(b) + \nu(b')$
and $\alpha(b + b') = \alpha(b) + \alpha(b')$. 
Therefore, $\psi$ is an isomorphism of $D$-algebras.
Finally,
denote by $E'_{A}$ the exponentiation on $D[t^{A}]$ with domain $D$
given by $E'_{A}(u + a) = E(u) t^{a}$ (for $u \in V$ and $a \in A$), and
similarly denote by $E'_{B}$ the exponentiation on $D[t^{B}]$ with domain $D$
given by $E'_{B}(u + b) = E(u) t^{b}$ (for $u \in V$ and $b \in B$).
Given $u \in V$ and $b \in B$, we have:
\begin{multline*}
\psi(E'_{B}(u + b)) = \psi(E(u) t^{b}) = E(u) E(\nu(b)) t^{\alpha(b)}
= \\
= E'_{A}(u + \nu(b) + \alpha(b)) = E'_{A}(u + b). 
\end{multline*}

2) Let $F=(F; V_F, E_F)$ be a partial E-domain extending $D = (D; V,E)$, and assume that $D\subseteq V_F$.
For every $a \in A$, define 
\[ \phi(t^{a}) \coloneqq E_{F}(a).\]
Extend $\phi$ to ${\bf D}'$ by $D$-linearity.
Then, $\phi$ is the unique homomorphism of E-domains over $D$ from ${\bf D}'$ into ~$F$.
\end{proof}

\medskip

Notice that a version of the above lemma holds also for a partial free 1-extension of $D,$ i.e.
$(D[t^{A}],V',E')$, where $V' = V \oplus A$ and A is not a complement of V inside D. 

\subsection{The free completion of a partial exponential ring}\label{sec:E-completion}

In this subsection, $R$ is a partial E-ring. 
\begin{definition}\label{def:RE}
For every $n \in \N$, define a chain of partial E-rings $(R_{n})_{n\in \mathbb N}$ inductively as follows
\begin{itemize}
\item 
 $R_{0} \coloneqq R$ 

\item 
$R_{n+1} \coloneqq  {\bf  R}_{n}'$ for all $n$
\end{itemize}

We define 
\[
R^{E} \coloneqq \bigcup_{n} R_{n}.
\]
 Notice that the exponential map on $R^{E}$ is
total.
The exponential ring $R^{E}$ is the free completion of~$R$.
The sequence $(R_{n})_{n \in \N}$ is the canonical filtration (over $R$) of
$R^{E}$ into partial exponential rings. 
\end{definition}
If the exponential function is total on $R$ then $R^{E}=R$.

\begin{remark}
$R^{E} = (R_{n})^{E}$ for every $n \in \N$.
\end{remark}

\begin{lemma}
Let $S$ be a partial E-ring.
Let $\phi: R \to S$ be a homomorphism of partial exponential rings.
Then, there exists a unique homomorphism of exponential rings
$\phi^{E} :R^{E} \to S^E$ extending $\phi$.
\end{lemma}
\begin{proof}
It is immediate by Lemma \ref{unico}.
\end{proof}

When the above map $\phi^{E} : R^{E} \to S^{E}$ is an isomorphism, we simply write $R^{E} = S^{E}$.

\begin{remark}
The rings $\C[\x]_{n}$ (for $n \ge 1$) and $\C[\x]^{E}$ are isomorphic 
\emph{as rings},
since they are both of the form $\C[\x][t^{G}]$ for some torsion-free divisible
Abelian group $G$ with cardinality of the continuum.  
\end{remark}

\bigskip
If $S$ is a subring of $R$ (or of $R^{E}$), we consider  $S$ as a partial
exponential subring in a natural way, by defining the domain of the
exponential map on $S$ as
\[
V_{S} = \set{a \in S: E(a) \in S}.
\]
In the following, unless  otherwise stated, we make the above canonical
choice for~$V_{S}$.

\bigskip

Given $R$ a partial E-ring and $S$ a subring of $R^E,$ we give some sufficient conditions on $S$ such that $S^E$ is isomorphic to $R^E.$
In order to do so we need the following notation.
Given $S$ a partial exponential subring of $R^{E}$, we write 
$S = R[e^{A}]$ if:
\begin{enumerate}
\item $R$ is a partial exponential subring of $S$; 
\item $A$ is a $\Q$-vector subspace of $R$;
\item $A \cap V = \{0\}$, where $V$ is the domain of the exponential map on $R$;
\item the canonical map $\theta: R[t^{A}] \to S$ sending $t^a$ to $E(a)$ for all $a\in A$ is an isomorphism of partial
E-rings.
\end{enumerate}
Notice that, when $S = R[e^{A}]$, we can identify $S$ with a partial exponential
subring of $R'$ (or, equivalently, of $R^{E}$); moreover, the domain of the exponential map on $S$ is
\[
V \oplus A = \set{x \in S: E'(x) \in S}.
\]


\begin{remark}
Let $S$ be a partial exponential ring such that $R \subseteq S \subseteq R^{E}$.
Then, there are two canonical maps $\alpha^{E}: R^{E} \to S^{E}$
and $\beta^{E}: S^{E} \to R^{E}$.
We have that $\beta^{E} \circ \alpha^{E} : R^{E} \to R^{E}$ is the identity map, and
therefore
$\alpha^{E}$ is injective and $\beta^{E}$ is surjective.
\end{remark}

We give now a sufficient condition guaranteeing that the above maps $\alpha^{E}$
and $\beta^{E}$ are isomorphisms (one inverse of the other), and therefore
$R^{E} = S^{E}$.

\begin{lemma}\label{lem:R-S}
Let $S$ be a partial exponential subring of  $ R^{E}$, and  
assume that $S = R[e^{A}]$ for some  subspace $A$ of the complement of $V$.
Then, $S^{E} = R^{E}$.
\end{lemma}
\begin{proof} First of all we prove the following 
\begin{claim}
$S' = R'[e^{C}]$ where $C$  is a $\Q$-vector subspace of $S$.
\end{claim}
Let $B < R$ be a $\Q$-vector subspace, such that
$V \oplus A \oplus B = R$.
Then,
\[
{R}' = R[t^{A \oplus B}] = R[t^{A}][t^{B}] = S[t^{B}].
\]
Remember that the domain of the exponential map on $S$ is
\[
V_{S} = V_{R} \oplus A.
\]
Let $C < S$ be a $\Q$-vector subspace, such that $R \oplus C = S$.
Then, 
\[
S = R \oplus C = V \oplus A \oplus B \oplus C = V_{S} \oplus B \oplus C.
\]
Therefore,
\[
S' = S[t^{B \oplus C}] = S[t^{B}][t^{C}] = R'[t^{C}].
\]

Using inductively the above claim, we obtain that, for every $n \in \N$,
$S_{n} = R[e^{C_{n}}]$ (for some suitable $C_{n}$). Therefore, the
canonical map $S^{E} \to R^{E}$ is injective, and hence it is an isomorphism.
\end{proof}

The following two examples show that in general the previous lemma is not true without the hypothesis on S. 

\begin{example}
Let $K$ be an exponential field, and $R \coloneqq K[x]$.

\begin{enumerate}
\item
Let $S_1 = R[e^{e^{x}}]$.
\item
Let $S_2 = R[e^{x} + e^{x^{2}}]$.
\end{enumerate}

As partial E-domains, $S_1$ and $S_2$ are both isomorphic to $K[x,y]$
(via the maps  $x \mapsto x$, $y \mapsto e^{e^{x}}$ and 
$x \mapsto x$, $y \mapsto e^{x} + e^{x^2}$, respectively), and therefore 
$S_i^{E}$ is isomorphic to $R[x,y]^{E}$, and hence $S_i^{E} \ne R^{E}$, 
for $i=1,2$. 
\end{example}

\begin{remark}\label{rem:standard-completion}
The construction of  the free completion $R^E$ of a partial E-ring $R$ is a
generalization of the well known construction of the ring of exponential
polynomials in any number of variables, see \cites{van, macintyre}. We
briefly recall it. 
Starting from an E-ring $R$, we construct  the E-polynomial ring in the
variables $\overline x=(x_1, \ldots, x_n)$, denoted by $R[\overline x]^E$.

Let $R_{-1} = R$, and $R_0 = R[\overline x]$ as partial exponential rings (where
$E$ is defined only on~$R$).
Then we can define $R[\overline x]^{E}$ as the free completion $R_{0}$ as in
Definition~\ref{def:RE}. 
In this case \cites{van,macintyre} give also a standard
way to choose the $\Q$-vector spaces $(A_{k}: k \in \N)$ such that, for every $k \in
\N$, $R_{k} = R_{k-1} \oplus A_{k}$,
$R_{k+1} = R_{k}[t^{A_{k}}]$, and $E(a) = t^{a}$ for every $a \in A_{k}$.
Let $A_0 = (\overline x) $ be the ideal
of $R[\overline x]$ generated by $\overline x$. 
For $k \geq 1$, let $A_{k}$
be the $R_{k-1}$-submodule of $R_{k}$  freely
generated by $t^a$ with  $0 \neq a \in A_{k-1}$.
Clearly,  $R_k = R_{k-1} \oplus A_k$, as $\Q$-vector spaces (and actually as
$R_{k-1}$-vector spaces) and the other conditions also hold. Sometimes it is convenient to represents $R[\overline x]^{E}$ as the group ring  $R[\overline x][t^{\bigoplus_{i\geq 0} A_i}]$. 
\end{remark}



In the construction we presented above we prefer to work with $D$ a $\mathbb
Q$-algebra, 
and any $\Q$-vector subspace as a
summand. 
Our construction gives also more flexibility in extending the domain of definition of the partial exponential function.

\section{E-ideals}\label{sec:E-ideal}
 We introduce the notion of E-ideals for partial E-rings. 
Let $D=(D; V,E)$ be a partial E-ring.

\begin{definition}
An E-ideal of $D$ is an ideal $I$ of the ring $D$ such that, for every
$v \in I \cap V$, $E(v) - 1 \in I$. 
\end{definition}

\begin{remark}
\label{trivialEideal}
If $I$ is an ideal of $D$ with $I\cap V=( 0)$ then $I$  is an E-ideal of $D$.
\end{remark}

In order to study E-ideals it is useful to recall the notion of augmentation map.

\begin{definition}
Let $D[t^A]$ be a free 1-extension of~$D$. The augmentation map $\sigma_{A}$ associated  to $A$ is the ``usual'' augmentation
map from $D[t^A]$ to $D,$  i.e. 
\[
\sigma_A(\sum d_i t^{a_i}) = \sum d_i .
\]
The corresponding augmentation ideal is 
\[
J_{A} := \sigma_A^{-1}(0) < D[t^{A}],
\]
\end{definition}
\noindent and it is an E-ideal. 

Notice that $\sigma_{A}$ and $J_{A} $  depend on the choice of $A$ (even when $A$
is a complement of $V$). In the sequel when clear from the context we  drop the subscript.

\begin{remark}
\begin{enumerate}
\item $J_{A} \cap D = \{0\}$.  
\item $J_{A}$ is an E-ideal of $D[t^{A}]$.  
\item As a
$D$-module, $J_{A} $ is generated by
\[
(E'(a) - 1 : a \in A)
\]
(where $E'(a) = t^{a}$).
\end{enumerate}
\end{remark}
\begin{proof}
(1) If $d \in D$ then $\sigma_{A}(d) = d$. 
Therefore, if $d \in J_{A} \cap D$, then $d = \sigma_{A}(d) = 0$.

(2) It follows immediately from (1).

(3) If $a \in A$, then $\sigma_{A} (t^{a} - 1) = 0$, and therefore $t^{a} - 1 \in J_{A} $.
Let $b \in J_{A}$, then $b$ can be written uniquely as 
 $b = \sum_{a \in A'} d_{a} t^{a}$, where $A'\subseteq A$, $|A'| < \infty$,  $d_{a} \in
D$, and $\sum_{a \in A} d_{a} = 0$.
Thus,
\[
b= \sum_{a \in A} d_{a} (t^{a} -1).
\]
\end{proof}

Also for (partial) E-rings there is a notion of prime and maximal  E-ideal.

\begin{definition} An E-ideal $I$ of $D$ is \intro{prime} if it is prime as ideal,
i.e.\ iff $D/I$ is a domain. An E-ideal $I$ of $D$ is an  \intro{E-maximal} ideal if it is maximal among the E-ideals. It is \intro{strongly maximal} if it is maximal as ideal.
\end{definition}


\begin{remark} Given an E-ideal $I$,  if $I$ is strongly maximal, then it is E-maximal and prime.  We show that the converse is not true. Moreover, 
neither E-maximal implies prime nor prime implies E-maximal.

Notice that on the quotient ring $D/I$ there is a natural partial exponential map E  whose domain of definition is $(V+I)/I$, and it is  defined as $E(v+I)=E(v)+I$. This is well defined.  Indeed, assume $v+I=v'+I$ then $v-v' \in I$, and since $I$ is an E-ideal we have $E(v-v')-1\in I$. So, $E(v)E(v')^{-1}-1=r$ for some $r\in I$, hence $E(v)-E(v')=rE(v')\in I$, and so $E(v)+I=E(v')+I$ (see also \cite{macintyre}). The  correspondence between ideals of $D$ containing $I$ and ideals of $D/I$ given by the canonical projection $\pi:D\rightarrow D/I$ extends easily to E-ideals. Moreover, $J$ is a maximal E-ideal of $D$ containing $I$ iff $\pi(J)$ is a maximal E-ideal of $D/I$. So, $I$ is E-maximal in $D$ iff $\pi(I)=0$ is E-maximal in $D/I$.

\end{remark}

Given an ideal $I$ of a ring $R$, there are two canonical ways to extend it to an ideal of the group ring
 $R[G]$, where $G$ is any Abelian group (see  also \cite{point}). Let $\sigma: R[G] \to R$ be the associated
augmentation map then for any ideal  $I$ of  $R$
we define the following two ideals of $R[G]$
\[
I'_{G} := \sigma^{-1}(I), \qquad
I''_{G} := I R[G],
\]
where  $ I R[G]$ is the ideal of $R[G]$ generated by $I$, i.e. 
\[ I R[G]=\{ \sum_{i=1}^m r_ig_i:  r_i\in I \mbox{ and } g_i\in G, m\geq 1 \}.\]

The ideals $I'_{G}$ and $I''_{G}$ satisfy the following properties.

\begin{lemma}\label{lem:I-ext}
\begin{enumerate}
\item
$I \subseteq I''_{G} \subseteq I'_{G}$. 
\item $I'_{G} \cap R = I=I''_{G} \cap R$.
\item
The map $\sigma$ induces an isomorphism between $R[G]/I'_{G}$ and
$R/I$.
\item
$I$ is prime (resp., maximal) in $R$   iff $I'_{G}$ is prime  (resp., maximal) in $R[G]$. 
\item If $G$ is torsion free and $I$ is prime, then $I''_{G}$ is prime.
\item $I'_{G}$ is maximal among the ideals $L$ of $R[G]$ such that
$L \cap R = I$.
\end{enumerate}
\end{lemma}
\begin{proof}
(1) It is trivial

\noindent (2) The two equalities  follow from  $\sigma(r) = r$ for every $r \in R$.


\noindent (3) The map from $R[G]$ to $R/I$ induced by $\sigma$ is a surjective ring
homomorphism whose kernel is $I'_{G}$, and therefore $R[G]/I'_{G}$ is
isomorphic to $R/I$.

\noindent (4) It is an immediate consequence of (3).

\noindent (5) If $I$ is prime then $R/I$ is an integral domain, and if moreover
$G$  is torsion-free, then  $(R/I)[G]$ is also an integral domain. The natural map
 \[
\begin{array}{ccccc}

\tau  & : & R[G]  & \rightarrow & (R/I)[G] \\

    &     &       &                           &                \\

     &    & \sum_{i=1}^m r_ig_i & \mapsto & \sum_{i=1}^m (r_i+I)g_i
  \end{array}
\]
 is a ring homomorphism  whose kernel is $I''_{G}$, and since $\tau$ is also
 surjective  the rings  $R[G]/I''_{G}$ and $(R/I)[G]$  are
 isomorphic. Therefore $I''_{G}$ is prime.

\noindent (6) 
Let $L$ be an ideal of $R[G]$ containing~$I'_{G}$. 
Suppose that there exists $b \in L \setminus I'_{G}$.
Then, $\sigma (b) \notin I$.
Moreover,
\[
\sigma (b) - b \in \sigma^{-1}(0) \subseteq I'_{G} \subset L.
\]
Therefore,
\[
\sigma(b) \in L \cap R \supsetneq I. \qedhere
\]
\end{proof}

\smallskip

We now consider an E-domain $D=(D; V,E)$.
Let $I$ be an ideal of $D$,   $A < D$ a complement of $V$ in $D$, and $V'=V\oplus A$.  Then 
 we can construct as above the two ideals 
$I'_{A}$ and $I''_{A}$ of  $D[t^{A}]$ extending $I$.
In general, they are not E-ideals.
However, in some specific cases if $I$ is an E-ideal, so also are $I'_{A}$ and $I''_{A}$. In the next lemmas we examine properties of $I$ which are inherited by $I_A'$ and $I_A''$. In particular, we analyze properties of  $I'_{A}$ and $I''_{A}$ for $I$ an E-ideal in two mutually excluding cases when $A\subseteq I$ and $(I+V)\cap A=\{ 0\}$. 
 
\begin{lemma}\label{lem:aug-extension}
\begin{enumerate}
\item $I'_{A} \cap D = I''_{A} \cap D = I$.
\item 
$I''_{A} \cap V' = I \cap V'=I'_{A} \cap V'$.
\item $(I_A' + V') \cap D = I + V+A$.
\item  $I$ is prime (resp. maximal) in $D$ if and only if  $I'_{A}$ is prime (resp. maximal) in $D[t^{A}]$.
\item If $I$ is prime, then $I''_A$ is prime.
\end{enumerate}
\end{lemma}

\begin{proof}
(1) This is clear from the definitions of $I'_{A}$ and $ I''_{A} $.

(2) This follows from (1) and the definitions of $I'_{A}$ and $ I''_{A} $.

(3) Let $b \in (I_A' + V') \cap D$.
Thus,  $b = b' + v'$, for some $b' \in I_A'$ and some $v \in V'$.
From  $V' = V \oplus A$, it follows that $v' = v + a$, for some $v \in V$ and $a \in A$. 
Moreover, $b$ and $v$ are in $D$, and this implies  $\sigma_{A}(b) = b$ and $\sigma_{A}(v) = v$.
On the other hand, since $a \in A$, $\sigma_{A}(a) = a$.
Finally, since $b' \in I_A'$, we have $\sigma_{A}(b') \in I$.
Therefore,
\[
b = \sigma_{A}(b) = \sigma_{A}(b') + \sigma_{A}(v) + \sigma_{A}(a)=  \sigma_{A}(b') + v +a\in I + V+A.
\]
The other inclusion is trivial.

(4) See (4) in Lemma~\ref{lem:I-ext}.

(5) See (5) in  Lemma~\ref{lem:I-ext} (recall that $A$ is torsion free).

%
%
\end{proof}


We now consider the case when $I$ is an E-ideal and $A\subseteq I$. 

\begin{lemma}\label{A contained in I}
 If $I$ is an E-ideal of $D$ containing $A$, then
$I'_{A}$ is an E-ideal of $D[t^{A}]$.
In this case $I'_{A}$ is the  E-ideal generated by $I$: in
fact, as a $D$-module, $I'_{A}$ is generated by $I$ and 
 $(E'(a) - 1: a \in A)$.
 
Moreover, if $I$ is E-maximal, then $I'_{A}$ is also E-maximal.

\end{lemma}

\begin{proof}
Let $r\in I'_{A}\cap V'$. Hence, 
$r=v+a$ for unique $v\in V$ and $a\in A$, and so $r\in D$. Since $r\in  I'_{A}$ we have  $\sigma_A(r)=r=v+a\in I$, and this implies $v\in I$, and so $E(v)-1\in I$. Now,  $\sigma_A(E'(r)-1)=\sigma_A( E(v)t^a-1)= E(v)-1\in I$, so $E'(r)-1\in I'_{A}$.

Let $b \in I'_{A}.$ Then we can write uniquely $b = \sum_{a \in A'} d_{a} t^{a}$ where $A'\subseteq A$, $|A'|<\infty$ for some $d_{a} \in D$, and $\sigma_{A}(b) \in I$.
Thus,
\[
b= \sum_{a \in A'} d_{a} (t^{a} -1) + \sigma_{A}(b).
\]
Therefore, $I'_{A}$, as a $D$-module, is generated by $I$ and 
 $(E'(a) - 1: a \in A)$ (recall that $E'(a)=t^a$ for all $a \in A$).



Assume, now,  that $I$ is E-maximal. We want to show that $I'_{A}$ is also
E-maximal.
Let $P$ be an E-ideal of $D[t^{A}]$ such that $P \supsetneq I'_{A}$.
By Lemma~\ref{lem:I-ext}, there exists
$b \in P \cap D \setminus I$. 
Since $I$ is E-maximal, $1$ is in the E-ideal generated by $I + b$, and
therefore
$1 \in P$, proving that $I'_{A}$ is E-maximal.
\end{proof}

In the next lemma  consider the case of $I$ an E-ideal and   $(I + V) \cap A = \{ 0\}$.
\begin{lemma}\label{lem:flat-extension}
Let $I < D$ be an E-ideal, and assume $(I + V) \cap A =\{ 0\}$.
Then
\begin{enumerate}
\item  $I''_{A} \cap V' = I \cap V$.
\item  $I''_{A}$ is an
 E-ideal, and  $I''_{A}$ is the E-ideal of $D[t^{A}]$
generated by~$I$. If $I$ is prime then $I''_{A}$ is prime. 
\item  $I'_{A}$ is also an E-ideal
containing $I$. 
If moreover,  $I$ is either prime, strongly maximal, or E-maximal, then
$I'_{A}$ is also prime, strongly maximal, or  E-maximal, respectively.
\end{enumerate}
\end{lemma}
\begin{proof}


(1) By (2) of Lemma \ref{lem:aug-extension} and the hypothesis  $(I + V) \cap A = \{ 0\}$,  we have
\[
I''_{A} \cap V' = I \cap V' = I \cap (V \oplus A) = I \cap V.
\]

(2)  Let $b \in I''_{A} \cap V'$.
By (1), $b \in I \cap V$, and  since $I$ is an E-ideal of $D$,
$E(b) -1 \in I \subseteq I''_{A}$. By (5) of Lemma \ref{lem:I-ext}, if $I$ is prime then $I''_A$ is prime. 

(3) We show that $I'_{A}$ is an E-ideal.
Let $b \in I'_{A} \cap V'$.
Then, $b\in V' = V \oplus A \subseteq D$ implies  $b = \sigma(b) \in I \cap (V \oplus A) = I \cap V$ (the last equality follows from the hypothesis $(I + V) \cap A = \{ 0\}$).
Since $I$ is an E-ideal of $D$, we have $E(b) - 1 \in I \subseteq I'_{A}$. 
The remaining claims are proved  in a similar way as  the corresponding claims
in  Lemma~\ref{lem:aug-extension} and Lemma~\ref{A contained in I}.

\end{proof}

In the proof of Lemma \ref{lem:R-S} we extended the domain of the exponential
function to $D$ into two steps according to a decomposition of $A$, the
complement of $V$ in~$D$. 
We now use this method in order to extend a fixed  E-ideal $I$ of $D$
to an E-ideal of the two (or more) free 1-extensions associated to a 
decomposition of~$A$.

Let $I < D$ be an E-ideal. Define $I_{0} := I \cap V$. Let $P_{1}$ be a complement of $I_{0}$ inside $I$ (as $\Q$-vector spaces)
and $P_{2}$ be a complement of $I +V$ inside $D$. We use the above complements to
give the following decomposition of~$A$, i.e.\  $A := P_{1} \oplus P_{2}$. 
Thus, we can build ${\bf D' }= D[t^{A}]$ the corresponding free 1-extension
of $D$, with exponential map~$E'$ in two steps.

Let $F := D[t^{P_{1}}]$ as partial E-domain, where the
corresponding exponential map is $E_{F}$, its domain is $W = V \oplus P_{1}$, and $E_F(v+b) =E(v)t^b$ for every $v\in V$ and every $b\in P_1$.
Since $A = P_{1} \oplus P_{2}$, as rings we have
$D' = F[t^{P_{2}}]$, and $E_{F}$ extends $E$ and is extended by~$E'$.  Let $\tau_1: D[t^{P_{1}}]\rightarrow D$ and $\tau_2: D[t^{P_{1}}][t^{P_{2}}]\rightarrow D[t^{P_{1}}]$ the corresponding augmentation maps. It is easy to verify that $\sigma_A=\tau_1\circ \tau_2$. \


\begin{lemma}
In the setting above: 
\begin{enumerate}
\item $W = V + I$ and $P_{2}$ is a complement of $W$ inside $D$. 
\item
\[
(I'_{P_{1}} + W) \cap P_{2} = 0
\]
\item $D' = F[t^{P_{2}}]$ as partial E-domains.
\end{enumerate}
\end{lemma}
\begin{proof}
(1) Since $V = V + I_{0}$, we have
\[
W = V  + P_{1} = V + I_{0} + P_{1}  = V + I.
\]
Moreover,
\[
P_{2} \cap W = P_{2} \cap (I + V) = 0,
\]
and 
\[
P_{2} + W =P_{2} + (I + V) = D.
\]
Thus, $P_{2}$ is a complement of $W$ inside $D$.

(2) Since $P_{2} \subseteq D$,  by  Lemma~\ref{lem:aug-extension}(1) we have $(I'_{P_{1}} + W) \cap P_{2} = (I + W) \cap P_{2} = \{ 0\}$.

(3) First, we show that the domain of the exponential function of $D'$ is the
same as the one of $F[t^{P_{2}}]$.
The former is $D$, while the latter is $V + P_{1} + P_{2}$, which is equal
to~$D$.
Second, we show that the two exponential functions coincide.
Let $a \in D$. 
We can write  $a = v + p_{1} + p_{2}$ for some unique $v \in V$,
$p_{1} \in P_{1}$ and $p_{2} \in P_{2}$.
We have
\[
E'(a) = E(v) t^{p_{1} + p_{2}} = E(v) t^{p_{1}} t^{p_{2}}
=E_{F}(v + p_{1}) t^{p_{2}}
\]
and the latter is equal to the exponential function in $F[t^{P_{2}}]$.
\end{proof}

We use the above decomposition to extend $I < D$ to an E-ideal of $D'$.
Define $H_{1} := \tau_{1}^{-1}(I)=I'_{P_{1}}$.
Since $P_1\subseteq I$, by Lemma~\ref{A contained in I}, $H_{1}$ is an E-ideal of~$F$.
Now, let $H := H_1 D'$, the ideal of $D'$ generated by $H_1$. Since $(H_1 + W)\cap P_2=\{ 0\}$,  Lemma~\ref{lem:flat-extension} implies that 
 $H$ is an E-ideal of~$D'$.

Let $L := \tau_{2}^{-1}(H_{1}) = \sigma_A^{-1}(I)$.
Combining Lemma~\ref{lem:aug-extension}
  and~\ref{lem:flat-extension}, we have the following lemma.

\begin{lemma}\label{lem:ideal-1-extension}
\begin{enumerate}
\item The ideal $H$ is the E-ideal of $D'$ generated by~$I$.

\item
$L = \sigma_A^{-1}(I)$ is the E-ideal of $D'$ generated by $I + \sigma_{A}^{-1}(0)$.

\item
If $I$ is a prime (as ideal), then $H$ and $L$ are prime.

\item
If $I$ is a proper ideal of $D$, then $H$ and $L$ are a proper ideals of~$D'$.

\item
If $I$ is strongly maximal (resp., E-maximal), then $L$ is strongly maximal
(resp., E-maximal).

\item
$H \cap D = L \cap D = I$.
\end{enumerate}
\end{lemma}

In general,  $H$ is not $ E$-maximal, even when $I$ is strongly maximal (it
suffices that $V + I \ne D$ to get that $H$ is not E-maximal).

If $V + I = D$ then $P_2=0$, i.e.\ $D'=D[t^{P_1}]$. 
So, $H=H_1$ which is strongly maximal, and so $H$ is strongly maximal by Lemma~\ref{lem:I-ext}.
  
Assume $V + I \neq D$. Then $I \leq H \le L$ 
and $P_2 \neq 0$;
let
$0 \neq a \in P_2$, then $t^a-1 \notin H$ otherwise $1\in H_1$, and $t^a-1\in L$, so $H \subsetneq L$, and so $H$ is not E-maximal. 

\medskip




By iterating Lemma~\ref{lem:ideal-1-extension} we obtain the E-ideal of $D^E$ generated by $I$ which we  denote by $I^{E}$. 
\begin{thm}\label{thm:ideal-complete-extension}
Let $I \subseteq D$ be an E-ideal.
Then
\begin{enumerate}
\item $I^E \cap D = I$;
\item $I$ is prime iff $I^E$ is prime;
\item if $I$ is E-maximal (resp., strongly maximal) then there exists $J'$ ideal of
$D^{E}$ which is E-maximal (resp., strongly maximal) and such that $J' \cap D = I.$
\end{enumerate}
\end{thm}

\section{Maximal and Prime E-ideals}

Clearly, a strongly maximal E-ideal is both prime and E-maximal;
in this section, we show independence of E-maximality, strong maximality, and
primeness of E-ideals, beyond the above trivial implication.
Our examples are constructed in the ring of exponential polynomials over an
E-field~$K$. 

Using Lemma~\ref{lem:R-S} we see that not every E-maximal ideal is prime. 

\begin{thm}
\label{maxnoprime}
There exists an E-ideal $P < K[x]^{E},$ such that $P$ is E-maximal but not
prime (and therefore $P$ is not strongly maximal).
\end{thm}

\begin{proof}
Let $S = K[x, e^{\Q x}]$. We consider $S$ as a partial exponential subring of
$K[x]^{E}$.
Let $I = (x^{2}, e^{x} + 1)$. Clearly,  $I$ is a proper ideal of the ring $S$.
Notice that $I$ is also a proper  E-ideal of $S$ (see Remark~\ref{trivialEideal}).
Let $I^{E}$ be the  E-ideal generated by $I$ inside $S^{E}$. 
By Theorem~\ref{thm:ideal-complete-extension}, 
$I^{E}$ is  a proper E-ideal of $S^{E}$.

Moreover, $S = K[x][e^{\Q x}]$, and therefore, by Lemma~\ref{lem:R-S},
$K[x]^{E} = S^{E}$, thus, we can identify $I^{E}$ with a proper E-ideal of $K[x]^{E}$
(the E-ideal generated by $(x^{2}, e^{x}+1)$ in $K[x]^{E}$).
Let $P < K[x]^{E}$ be an E-ideal of $K[x]^{E}$ containing $(x^{2}, e^{x}+1)$, and E-maximal.
We claim that $P$ is not prime.
Indeed, if $P$ were prime, then $x \in P$, and therefore $e^{x}-1 \in P$, and
therefore $2 \in P$, contradicting the fact that $P$ is a proper ideal.
\end{proof}


We can also produce an E-ideal which is prime, E-maximal, but not strongly maximal.

\begin{thm}\label{lem:emax-not-max}
There exists a (proper) E-ideal $P$ of $K[x]^{E}$ such that $P$ is prime,
E-maximal, but not strongly maximal.
\end{thm}

\begin{proof}
Let $R = K[x]$.
Let $A := xR$ (the ideal generated by $x$), 
and $B$ be a $\Q$-linear basis of $A$, such that
$x, x^{2}, x^{3}, \dots \in B$.
Let $U := R \setminus A$, and $\alpha: B \to U$ be some 
surjective map
such that
$\alpha (x^{n}) = x^{n} -1 $ for every $1 \leq n \in \N$ (the reason of this choice will be clear at the end of the proof).
Let $F$ be the fraction field of~$R$, 
$\bar F$ be its algebraic closure, and $\bar F^{*}$ be the group of nonzero
elements of $\bar F$.
Let $\beta: (A, +) \to (\bar F^{*}, \cdot)$ be some group homomorphism extending $\alpha$.
By the universal property of group rings, 
there exists a unique  homomorphism of $R$-algebras 
\[\gamma: R_{1} = R[t^{A}] \to \bar F
\]
such that, for every $a \in A$, $\gamma(t^{a}) = \beta(a)$, and $\gamma\vert_{R}=id_R$. 
Let $P_{1} := \ker(\gamma) \subseteq R_{1}$, and  $S := \gamma(R_{1}) \subseteq \bar F,$ which is isomorphic to the quotient $R_1 / P_1$.
It is enough to prove that  $P_{1}$ is an E-ideal of $R_{1}$
which is prime, E-maximal, but not strongly maximal. Indeed, 
by Theorem~\ref{thm:ideal-complete-extension}, 
$P_{1}$ can then be extended to an E-ideal $P$ of $R^{E}$ 
satisfying the conclusion.

Since $\gamma(1) = 1$,  $\gamma$ is not identically zero, and thus $P_{1}$ is a proper ideal of $R_{1}.$
Moreover, since $\gamma$ restricted to $R$ is an embedding then $P_{1} \cap R = (0)$.
By Remark~\ref{trivialEideal} we have that  $P_{1}$
is an E-ideal.
Moreover, $S$ is a domain since $\bar F$ is a field,  and therefore
$P_{1}$ is prime.

It remains to show that $P_{1}$ is  not strongly maximal but E-maximal.
Thus, we have to show that $S$ is not a field and $(0)$ is an E-maximal ideal
of $S$.\\
In order to prove that $S$ is not a field we show that there are nonzero elements in $S$ which are not invertible.  Let $T := R_{(x)} = R[U^{-1}]$ be the localization of $R$ at the ideal $xR$.
Notice that $S$ is contained in the integral closure of $T$ (since it is
obtained by adding  $n^{th}$ roots of elements of~$T$). Moreover, $T$ is a principal ideal domain and therefore it is a normal domain
(i.e.\ integrally closed in its fraction field, see \cite{Matsumura}). Therefore,
since $x$ is not a unit of $T$, $x$ is also not a unit of~$S$, hence $S$ is not a field.
This proves that $P_1$ is not strongly maximal.\\
Now we want to prove that $(0)$ is the only proper E-ideal of $S.$ 
Let $J \leq S$ be some non-zero E-ideal. We have to show that $J = S$.
Let $0 \neq a \in J$, and let $p(z) \in T[z]$ be a nonzero monic polynomial of minimum
degree such that $p(a) = 0$.
Write $p(z) = z q(z) - t_{0}$ for some $q \in T[z]$ and $t_{0} \in T$.
Notice that $t_{0} \neq 0$ (since $p$ is of minimum degree).
Moreover, since $p(a) = 0$, we have 
\[
t_{0} = - a q(a) \in a S \subseteq J.
\]
Thus, $0 \neq t_{0} \in J \cap T$.
Therefore, $J \cap T$ is a non-zero ideal of $T$. All the non-zero ideals of
$T$ are of the form $x^{n} T$ for some $n \in N$, and thus $x^{n} \in J$ for some $n
\in \N$.

If $n = 0$, then $1 \in J$, and therefore $J = S$.

If $n \geq 1$, then $E(x^{n}) - 1 \in J$ (since $J$ is an E-ideal).
By definition of $P_{1}$,
\[
E(x^{n}) - 1 =  (x^{n} - 1) - 1 = x^{n} - 2.
\]
Thus, $2 \in J$, and therefore $J = S$. This means that $(0)$ is the  E-maximal ideal of $S$ and so $P_1$ is E-maximal in $R_1.$ In fact suppose that there exists a proper E-ideal $J_1$ of $R_1$ which contains $P_1.$ From the fact that $(0)$ is the only proper E-ideal of S it follows that $P_1 = J_1$ which conclude the proof.

\end{proof}

\section{Unique Factorization domain }\label{Hilbert_Null}

In this section we analyze classical algebraic properties  satisfied by polynomial rings over a field in the context of the ring of exponential polynomials over an E-field. 
So, we assume always  that $K$ is an E-field.

First of all we recall some notions and known results. An exponential polynomial $f(\x) \in K[\x]^E$ is invertible if it is of the form $f(\x)=e^{g(\x) }$ for some $g(\x) \in K[\x]^E$.  The notions of prime and irreducible element in $K[\x]^E$ are the usual ones.


 The pioneering work of Ritt (see \cite{Ritt}) on the factorization of the ring of exponential polynomials of height $1$ (see below for the definition of height) over $\C$ has been extended over the years to larger classes of exponential polynomial rings (see \cite{HRS}, \cite{GP74}), and more in general to group rings of a divisible  torsion free Abelian group  over a unique factorization domain (see \cite{Everest_Porten}).  The following factorization result follows from Theorem 1 of  \cite{Everest_Porten}) once the E-ring $K[\x]^E$ is considered as the group ring $K[\overline x][t^{\bigoplus_{i\geq 0} A_i}]$, see Remark~\ref{rem:standard-completion}.
 

\begin{thm}
Let $(K,E)$ be an E-field of characteristic $0$ with a non trivial exponential function E. Let $f(\x) \in K[\x]^E$, and $f \neq 0.$ Then $f$ factors, uniquely up to units and associates, as a finite product of the form: $$f = E(u) \cdot F_{1} \cdot \ldots F_{h} \cdot G_{1}(E(b_1)) \cdot \ldots \cdot G_{k}(E(b_k)),$$ where $F_{1}, \ldots, F_{h}$ are irreducible exponential polynomials  in $K[\x]^E$, and  $G_{1}, \ldots, G_{k}$ are non constant polynomials in one variable over $K$ and $u, b_1, \ldots, b_k \in A_{i} $, with the property that no $b_i$ can be written as a rational multiple of any $b_j$ with $i \not = j.$
\end{thm}


\begin{lemma}\label{irriducibilieprimi}
The ring of exponential polynomials $K[\overline x]^E$ is a GCD-domain. Therefore any irreducible element of $K[\overline x]^E$ is prime,  and so generates a prime ideal. 
\end{lemma}

\begin{proof} 
Since $K[\overline x]^E$ is a group ring over a divisible  torsion free Abelian groups, then by \cite[Theorem~5.2]{GP74} it is a GCD-domain. Hence, every irreducible element is prime. 
\end{proof}


Following the notation introduced in Remark~\ref{rem:standard-completion} to any $f \in K[\overline x]^E$ we associate its {\bf{height}}, denoted by $ht(f)$, i.e.  the smallest $n \in
\N$ such that $f \in R_{n}$. 

We say  that $f$ is {\bf reduced} if every associate of $f$ has height at least
$ht(f)$.

An example of an exponential polynomial which is irreducible and reduced is $x$;
an example of an exponential polynomial which is irreducible but not reduced is
$x e^{x}$, since the latter is an associate of $x$ and $ht(x) < ht(x e^{x})$.

\begin{lemma}\label{lem:reduced}
 Let $f, p$ be in $K[\x]^{E}$. 
If $f$ is a reduced polynomial  then $ht(fp) \geq ht(f)$.
\end{lemma}

\begin{proof} Suppose $f, p \in R_{n} \setminus R_{n-1}$.
Recall that  $R_{n} = R_{n-1}[t^{A_{n-1}}]$, 
where $A_{n-1}$ is a torsion free ordered Abelian group. So we can write uniquely 

$$f = q_1 e^{g_1} + \ldots + q_k e^{g_k} \mbox{ where } q_i \in R_{n-1} \mbox{ and } g_i \in {A_{n-1}},$$ and 
$$p = s_1 e^{b_1} + \dotsb + s_l e^{b_l}  \mbox{ where } s_i \in R_{n-1} \mbox{ and } b_j \in {A_{n-1}}.$$

Let $<$ be the order defined on  $A_{n-1}$. We can assume that

$$g_1< g_2 < \dotsb < g_k,$$ 

and

$$b_1< b_2 < \dotsb < b_l.$$

So we have $fp = \sum_{i =1}^{k} q_i(\sum _{j =1}^{l}  s_j e^{g_i + b_j}).$ If $ht(fp)
< ht(f)$, since we have $g_1 + b_1 < g_k + b_l$ then necessarily $g_i + b_j = 0$
for all $i = 1, \dotsc  k$ and $j = 1, \dotsc, l.$ This implies $g_1 = \ldots = g_k$
and $b_1 = \ldots = b_l.$ 
So, $f = (q_{1} + \dotsb + q_{k}) e^{g_{1}}$ is not reduced, because
$q_{1} + \dotsb + q_{k}$ is an associate of $f$ of strictly smaller height,
 in contradiction with the hypothesis.
\end{proof}

In the following result we prove that irreducible exponential polynomials generate not only ideals which are prime but also exponential ideals which are prime.

\begin{proposition}
\label{E-prime}
Let $f \in K[\x]^{E}$ be irreducible.  Let $(f)^{E}$ be the E-ideal of $ K[\x]^{E}$
generated by $f$.  Then, $(f)^{E}$ is a prime ideal in $K[\x]^{E}$. 
\end{proposition}
\begin{proof}
Let $n$ be the height of $f$ (that is, $f \in R_{n} \setminus R_{n-1}$).  
Without loss of generality, $f$ is reduced.

Then,
$f R_{n}$, the ideal of $R_{n}$ generated by $f$, is an E-ideal (see Remark~\ref{trivialEideal}), and it
is prime in $R_{n}$, by Lemma \ref{irriducibilieprimi}.
By (2) of Theorem~\ref{thm:ideal-complete-extension}, 
$(f)^{E}$ is a prime ideal in $K[\x]^{E}$.
\end{proof}

\noindent{\textsc{Open  question:}} When an E-prime ideal is finitely generated?



\section{Hilbert rings and Nullstellensatz}
Recall the following notion:

\begin{definition} An Hilbert ring (also Jacobian ring) is a ring where any prime ideal is the intersection of maximal ideals.
\end{definition}

In \cite{point}  the authors prove that:

\begin{proposition} The ring $R[\overline x]^E$ over any E-ring $R$ is not a Hilbert ring.
\end{proposition}

In the next proposition we show that the relativized property of being a {\it Hilbert ring} to  E-ideals is also not true for $K[\bar x]^{E}$. 


\begin{proposition}\label{lem:notHilbert}
Let $(K, E)$ be an exponential field such that E is non-trivial.
Then, there exists a prime E-ideal $P$ of $K[\bar x]^{E}$ such that
$P$ is not E-maximal and it is contained in exactly one E-maximal E-ideal.
In particular, $P$ is not an intersection of E-maximal E-ideals. 
\end{proposition}
\begin{proof}
For the sake of simplicity  we give the proof in the case  $\bar x=x$, noticing that the proof for the general case can be easily obtained from it.  

Let $\sigma: K[x]^{E} \to K[x]$ be the augmentation map associated to the E-ring $K[ x]^{E}$ as  group ring on the ring $K[ x]$,
and let $P$ be the augmentation ideal, i.e. $P:= \sigma^{-1}(0)$. Without loss of generality we can assume that an E-polynomial in $P$ has the following expression $$p_1(x)e^{xg_1(x)}+ \ldots + p_r(x)e^{xg_r(x),}$$ where 
$p_1(x)+ \ldots + p_r(x)=0$. 
Since $(0)$ is a prime  E-ideal in $K[x]$, $P$ is a prime E-ideal in
$K[x]^{E}$.

We prove that $I_{0} = \set{f \in K[x]^{E}: f(0) = 0}$ is the only proper
E-ideal properly containing $P$.  It is clear that $I_{0} \supseteq P$. Moreover, $P \neq I_{0}$, since $x\in I_0$  and  $x\not\in P$.

Let $J$ be any proper E-ideal such that $J$ properly contains~$P$. It is enough to show that if  $f\in J \setminus P$ then $f\in I_0$. This implies that $P$ is contained in no E-maximal ideal different from $I_0$.

We claim that it is enough to show that for any $f\in  J \setminus P$  and $f\in  K[x]$ then $f\in I_0$. Indeed, 
assume  $f$ is any element of $ K[x]^{E}$,  $ f\in J \setminus P$ and let $h=\sigma (f)$. Then $\sigma (f-h)= \sigma (f)-\sigma (\sigma (f))= h-h=0,$
so $g=f-h \in P$.  Moreover, $f - g=h \in K[x]$, and $f-g\in J\setminus P$. If $h\in I_0$ then also $f\in I_0$.

Let $p \in K[x] \cap (J \setminus P)$, and we prove that $p(0) = 0,$ i.e. $p \in I_0.$
Assume by contradiction that $c := p(0)\ne 0$. 
Since $E$ is not the trivial exponentiation taking constant value  $1$, there exists $0 \neq d \in K$ such that $E(c d)
\neq 1 $.
By replacing $p$ with $d p$, without loss of generality we may assume that $E(c) \neq 1$.
Let $p = c + q$, with $q \in (x)$.
By definition of $P$, we have 
\[
E(q) - 1 \in P \subset J
\]
and therefore
\begin{equation}\label{eq:cq-1}
E(c)E(q) - E(c) \in J;
\end{equation}
moreover, since $J$ is an E-ideal and $c + q = p\in J$,
\begin{equation}\label{eq:cq-2}
E(c) E(q) - 1 = E(p) - 1 \in J.
\end{equation}
Combining \eqref{eq:cq-1} and \eqref{eq:cq-2} we have
\[
E(c) - 1 \in J.
\]
Thus, $E(c) - 1$ is a non-zero element of the field $K$ which is in~$J$, and
therefore $J$ cannot be a proper ideal.


\end{proof}

In general, Hilbert Nullstellensatz is not true in the ring of exponential polynomials $K[\bar x]^{E}$, as shown in the following result.

\begin{proposition}\label{Null}
The exponential ring $\C[x]^E$  does not satisfy the weak Nullstellensatz, i.e.\ 
there exist two exponential polynomials $f$ and $g$ with no zeros in common, such
that the E-ideal generated by $f$ and $g$ is non-trivial.
\end{proposition}

\begin{proof} 

Let $f(x) = e^{x} - 1$ and  $g(x) =  e^{a x - b} - 1 \in \C[x]^E,$  where 
$a \notin \Q$ and $b \notin 2 \pi i (\Z + a \Z)$.

Let $J$ be the  E-ideal generated by $f$ and $g$ inside $\C[x]^{E}$.
Notice that $f$ and $g$ have no common zeros.
Indeed, if $f(c) = g(c) = 0$, then $c \in 2 \pi i \Z$ and $ac - b \in 2 \pi i \Z$,
contradicting the choice of~$b$.

It remains to prove that $J \neq \C[x]^{E}$.
Let $S := \C[x][e^{\C x}] \subset \C[x]^{E}$.
By Lemma~\ref{lem:R-S}, $\C[x]^{E}$ and $S^{E}$ are canonically isomorphic.
Thus, it suffices to show that $J \neq S^{E}$.
Let $L \subseteq S$ be the E-ideal of $S$ generated by $f, g$. Notice that $L$ is an E-ideal of $S$, and $J=L^E$.
It suffices to show that $L$ is a proper E-ideal of $S$ since by
Theorem~\ref{thm:ideal-complete-extension}, 
$J \cap S = L$, and we are done.

So, it suffices to prove that
$L \cap \C[x] = (0).$
Let $(a_{i})_{i \in I}$ be a $\Q$-basis of $\C$, with $a_{0} = 1$ and $a_{1} = a$.
As a $\C$-algebra, $S$ is isomorphic to $\C[x][y_{i}^{\Q}: i \in I]$ via an
isomorphism $\phi$ mapping $e^{a_{i} x }$ to $y_{i}$.
Let $r \in L \cap \C[x]$.
Thus, 
\[
r = p f + q g
\]
for some $p,q \in S$.
Let $R, P, F, Q, G \in \C[x][y_{i}^{\Q}: i \in I]$ be the images of $r, p, f, q, g$
respectively, under the map~$\phi$.
Notice that $R= r$ (so $R$ depends only on $x$), $F = y_{0} - 1$, $G= y_{1}e^{-b} - 1$.
So,
\[
r = P(x, \y) \cdot (y_{0}- 1) + Q(x, \y) \cdot (y_{1} e^{-b} - 1)),
\]
which holds in $\C[x][y_{i}^{\Q}: i \in I]$. 
If we specialize the above at $y_{0} = 1$ and $y_{1} = e^{b}$, we obtain $r = 0$, since $r$ does not depends on $\y$ but only $x$.
\end{proof}

\begin{remark}
The previous result holds also in a more general context:  for example,  the exponential polynomial ring over a Zilber field does not satisfy the weak Nullstellensatz.
\end{remark}

\section{Noetherianity}

As Macintyre points out in  \cite{macintyre},  
neither $\C[x]^E$ nor $\R[x]^E$ is Noetherian for E-ideals.  Indeed, the E-ideal 
$I = (E(\frac{x}{2^n}) - 1)_{n \in \N} $ 
is not finitely generated (for details see also \cite{terzo}).

\smallskip

There is a notion of noetherianity also for topological spaces.
\begin{definition}
A  topological space $X$ is Noetherian if it satisfies the descending chain condition for closed subsets, \ie any strictly descending sequence of closed subsets  of $X$ is stationary. \end{definition}

The E-ring $\C[\x]^{E}$ is very far from being Noetherian for E-ideals, since
$\C$ with the exponential Zariski topology is not a Noetherian space since 
the chain of closed subsets 
\[
\begin{array}{ccccccccc}
Z(E(x) -1) & \supset  & Z(E(\frac{x}{2}) -1) & \supset & \ldots  & \supset  & Z(E(\frac{x}{n!}) -1) &\supset  & \ldots\\

   &   &  &            &                &               &              &         &  \\
2\pi i\, \Z &\supset & 4\pi i\, \Z & \supset & \ldots & \supset & 2 n!\, \pi i\, \Z & \supset  &\ldots 
\end{array}
\]
descends indefinitely.
 
\smallskip

We show that $\C [\x]^E$ is not Noetherian even with respect to prime E-ideals.

\begin{thm}
The ring $\C [\x]^E$ does not satisfy the ascending chain  condition for prime E-ideals.

\end{thm}

\proof
We consider the subring $S := \C[\x, e^{\C \x}] =  \C[\x] [e^{\C\x}]  $
of the ring $\C [\x]^E$. 
The idea is to construct an ascending chain of prime
E-ideals of $S$ and to extend it to be an ascending chain ($P_i$) of prime E-ideals of  
$\C [\x]^E,$
by defining $Q_i = P_i  \C [\x]^E.$

For simplicity, we consider the case when we have only one variable 
$\x = x$.
We construct prime ideals in the following way. Let $B= (b_j : j < 2^{\aleph_0})$ be
a transcendence basis of $\C$ over $\Q$. For all $i \in \N,$ we define 
\[ p_i := e^{b_{3i}x} +e^{b_{3i + 1}x} + e^{b_{3i +2}x}  .
\]
Let  $A_n = (p_0, \ldots, p_n)$ be the ideal $\C[e^{\C x}]$ generated by $p_0, \ldots, p_n$.
We need the following 

\smallskip

\noindent 
{\bf Claim. }
The ideal $A_n$ is prime for all $n \in \N$.

\smallskip

\noindent {\bf Proof of the Claim.}
We introduce  new variables $z_i$ denoting elements of the form $e^{b_i x}$ for any
$b_i \in B$, \ie  $z_i = e^{b_i x}$. So,
$A_n = (z_0 + z_1 + z_2, z_3 + z_4 + z_5, \ldots, z_{3n} + z_{3n+1} + z_{3n + 2})$ and it is 
an ideal of $\C [\overline z^{\Q}].$
We prove by induction on $n$ that $A_n$ is prime. For $n=1$ we have that $A_1 =
(z_0+z_1+z_2)$ is an ideal of $\C [z_0^{\Q}, z_1^{\Q}, z_2^{\Q}]$. 
Assume that $p(z_0, z_1, z_2)\cdot q(z_0, z_1, z_2) \in A_1$ then $p(z_0, z_1, z_2)\cdot
q(z_0, z_1, z_2) = r(z_0, z_1, z_2) (z_0+z_1+z_2)$. Let $k$ be the common
denominator of any exponents in $p, q, r$, so we can consider $p, q, r \in \C
[z_0^{\pm \frac{1}{k}}, z_1^{\pm \frac{1}{k}}, z_2^{\pm \frac{1}{k}}]$. 
By replacing $z_{i}$ with $t_{i}^{k}$, we have that
$p,q,r \in \C[t_{0}^{\pm 1}, t_{1}^{\pm 1}, t_{2}^{\pm 1}]$.
Notice that $z_{0} + z_{1} + z_{2}$ becomes $s(\tv) := t_{1}^{k} + t_{2}^{ k} +
t_{3}^{k}$, and that $s$ is irreducible in $\C[\tv]$. Hence the ideal
$A_{1}'$ generated by $s$ inside $\C[\tv]$ is prime.
By localization, the ideal $A_{1}''$ generated by $s$ inside $\C[\tv^{\Z}]$ is
also prime; therefore, since $p q \in A_{1}''$, we have either $p \in A_{1}'' \subseteq
A_{1}$ or $q \in A_{1}'' \subseteq A_{1}$, proving that 
$A_1$ is prime.

For $n = 2$, $A_2 = (z_0+z_1+z_2, z_3+z_4+z_5)$; it is prime since 
\[
\frac{\C [z_0^{\Q}, \ldots, z_5^{\Q}]}{A_2} \cong \frac{\C [z_0^{\Q}, z_1^{\Q}, z_2^{\Q}]}{A_1} \otimes
\frac{\C [z_3^{\Q}, z_4^{\Q}, z_5^{\Q}]}{(p_1)}.
\]

So $\frac{\C [z_0^{\Q}, \ldots, z_5^{\Q}]}{A_2}$ is a domain, since the
tensor product of $\C$-algebras which are domains is also a domain, see
\cite[Chapter V, \S17]{bourbaki}. 
In a similar way we can prove that $A_n$ is prime for any $n \in \N.$ 
\qed

\smallskip

The Claim implies that   $P_n = A_n S$ is a partial prime E-ideal in S for each $n\in \N$. By Lemma~\ref{lem:R-S} and Theorem~\ref{thm:ideal-complete-extension} we have that also $Q_n$ is a prime E-ideal in  $\C [\x]^E$ for all $n\in \N$.  So we have a non stationary ascending chain of prime E-ideals.
\qed

\section{Strongly Maximal E-ideals \texorpdfstring{$I_{\overline a}$}{Ia}} \label{sec:smideal}
In this section  we investigate strongly maximal E-ideals. Let $(K,E)$ be an exponential field,
 $\overline a \in K^n$, and  $$I_{\overline a} = \{ f(\overline x) \in K[\overline x]^E : f(\overline a) = 0\}.$$

\begin{lemma}
\label{maximal}
  $I_{\overline a} $ is an E-ideal of $K[\overline x]^E.$  Moreover, it is a strongly maximal E-ideal.
\end{lemma}

\begin{proof}
It is immediate that  $I_{\overline a} $ is an E-ideal, and moreover it is
prime. Suppose, now,  $I_{\overline a}\subset J\subseteq K[\overline x]^E$ and let
$f(\overline x) \in J \setminus I_{\overline a}$, i.e. $f(\overline a)=\alpha\not=0$. 
Then $g(\overline x) := f(\overline x)-\alpha\in I_{\overline a}$, 
hence $\alpha\in J$, and so $J=K[\overline x]^E$.

\end{proof}

\begin{proposition} For all $\overline a \in K^n$,   $I_{\overline a}$ is $(x_1 - a_1, \ldots ,x_n - a_n)^E.$
\end{proposition}

\begin{proof}
Let $f(\overline x) \in I_{\overline a}$, we want to show that $f(\overline x) \in (x_1 - a_1, \ldots ,x_n - a_n)^E.$ 
We proceed by induction on the height of the polynomial. If the height of $f(\overline x)$ is zero, this means that it is a classical polynomial, and so $f(\overline x)\in  (x_1 - a_1, \ldots ,x_n - a_n)^E.$ Suppose  the result is true for all exponential polynomials of height at most  $k-1$, we prove the claim for $f(\overline x)\in   K[\overline x]^E $  where  height of $f$ is $k$. If $ht(f) = k$ then $f(\overline x) \in R_k  \setminus  R_{k-1}.$ Recall that  $R_k = R_{k-1}[t^{A_{k-1}}],$  and so $$f(\overline x)=\alpha_1 (\overline x)e^{\beta_1 (\overline x)} +\ldots + \alpha_m (\overline x)e^{\beta_m (\overline x)},$$ 
where $ht(\beta_i)<ht(f)$ for all $i=1, \ldots , m$. Suppose $f(\overline 0)=0$, and let $b_i=\beta_i (\overline 0)$, for $i=1,\ldots, m$. Let $\beta'_i (\overline x)= \beta_i (\overline x)-b_i$. Then $\beta'_i (\overline 0)=0$, and by inductive hypothesis $\beta'_i (\overline x)\in (\overline x)^E$. The polynomial $f(\overline x)$ can be written as follows, $$ f(\overline x)=\alpha'_1 (\overline x)e^{\beta'_1 (\overline x)} +\ldots + \alpha'_k (\overline x)e^{\beta'_k (\overline x)}$$ where $\alpha'_i (\overline x) = \alpha_i (\overline x)e^{b_i}$. Let $z_i=e^{\beta'_i (\overline x)} $, hence $f(\overline x)=g(\overline x,\overline z)$, and $g(\overline 0,\overline 1)=0$. By inductive hypothesis (ht($g$) < ht($f$)), $g(\overline x,\overline z)\in (\overline x,\overline z-\overline 1)^E$. From $\beta'_i (\overline x)\in (\overline x)^E$ it follows $e^{\beta'_i (\overline x)}-1\in (\overline x)^E$. So, 
$f(\overline x)\in (\overline x)^E$ ($z_i= e^{\beta'_i (\overline x)})$. A change of variables implies that if $f(\bar a)=0$ then $f\in (\x-\bar a)^E$.
\end{proof}

\begin{proposition} 
Let  $K$ be an algebraically closed E-field.
Any  prime E-ideal I of $K[x]^E$ containing a classical polynomial $p(x)$ is of the form $(x - a)^E,$ where $a$ is  a root of $p(x).$ Hence,  I is a strongly maximal E-ideal. 
\end{proposition}

\begin{proof}
Since $K$ is algebraically closed there are $\alpha_1,\ldots , \alpha_n\in K$ such that $p(x) =(x-\alpha_1)\cdot \ldots \cdot (x- \alpha_n)$. Hence, $x-\alpha_i\in I$ for some $i=1,\ldots, n$, and so $(x-\alpha_i)^E\subseteq I$. Lemma \ref{maximal} implies $(x-\alpha_i)^E= I$.
\end{proof}

\begin{remark}
The notions of prime E-ideal and E-maximal ideal are independent by Proposition
\ref{maxnoprime} and the following example. Let $q(x)\in K[x]^E$ be irreducible,
of height at least $1$ and reduced. By Lemma \ref{E-prime}, $(q)^E$ is
prime. Let $a\in K$ be a root of $q(x)$ (remember that e.g. on $\C$ all such
exponential polynomials $q$ have at least one root: see \cite{HR};
Example~\ref{ex:prime-not-max} gives an instance).

Then $q\in (x-a)^E$, and so $(q)^E \subset (x-a)^E$ (notice the proper inclusion since $x-a \not\in (q)^E$ because of the heights). So, $(q)^E$ is not E-maximal.
\end{remark}

\begin{example}\label{ex:prime-not-max}
Let $f = e^{\pi x} + e^{(\pi + \pi^2)x} + e^{\pi^2 x} -3 \in \C[x]^E$.  It  is irreducible in $\C[x]^E$.  Therefore $(f)^E$ is prime, and it is properly contained in $(x)^E$.
\end{example}

\section{Construction of maximal E-ideals}\label{sec:smideals-counterexample}

We now show that not all strongly maximal E-ideals are of the form $I_{\overline a}$ for some $\overline a \in K.$ We refer to the notation in Remark \ref{rem:standard-completion}, and 
we recall that $K[\overline x]^E = \bigcup_n R_n$.



\begin{lemma} \label{ideali}
Let $k\geq 1$. If $J_k$ is a maximal ideal of $R_k$ and $J_k \cap R_{k-1} = (0)$ then there exists a strongly maximal E-ideal $J$ of $K[\overline x]^E$ such that $J_k =J\cap R_k$. 
\end{lemma}

\begin{proof}
By Remark  \ref{trivialEideal},  $J_k$ is an E-ideal. For $n\geq 1$ each $R_n$ is a group ring over $R_{n-1}$ with associated augmentation maps $\sigma_{n}:  R_{n} \rightarrow R_{n-1}$. 


Let $\sigma_{k+1}$ be the augmentation map $\sigma_{k+1}:  R_{k+1} \rightarrow R_{k}$. 
Define  $J_{k+1} = \sigma_{k+1}^{-1}(J_k)$; 
 by  Lemma~\ref{lem:I-ext},  $J_{k+1} $  is a maximal ideal of $R_{k+1}$. Moreover, $J_k \subseteq J_{k+1}$ and $J_{k+1} \cap R_k = J_k,$ so $J_{k+1}$ is an E-ideal of $R_{k+1}.$ Iterating this construction we  define  $J = \bigcup_{n\geq k} J_n$ which is trivially an E-ideal of $K[\overline x]^E.$ Moreover, it is a maximal ideal, since if $a \not \in J,$ then $a \not \in J_n$ for all $n \geq k.$ Choose $n\geq k$ such that $a\in R_n$;  since $J_n$ is maximal then $J_n + (a) = R_n,$ and this implies that $J + (a) = K[\overline x]^E,$ i.e.\ $J $ is maximal.
\end{proof}

\begin{proposition}\label{thm:max-nontrivial}
 Let $K$ be an algebraically closed E-field with $|K| > \omega.$ Then there exists a
 strongly maximal E-ideal $J$ of $K[ x]^E$  such that $J \neq I_a$  for all $a \in K$.
\end{proposition}

\begin{proof}
Let  $S =K[x, E(cx) : c \in K]$, and $B$ a $\Q$-linear base of~$K$. 
Since
$\card K > \omega$ then $\card B = \card K$. 
We introduce new variables denoting elements of $E(Bx)$, one for each element
of $K$, i.e.\  $E(Bx) = \{y_c : c \in K\}$. 
Notice that all elements of $E(Bx)$ are invertible in~$S$. 
So,   $S = K[x, y^{\Q}_{c} : c \in K] $ (as a K-algebra). 
Clearly, $K[x] \subseteq S$. 
Consider the ideal $J _1= (y_c - (x-c) : c \in K)$ of $S$.   
Let $F := S/ J_1 = K[x, (x-c)^{\Q} : c \in K]$ be the quotient.
Since $K$ is algebraically closed, $K(x)$ is a field; moreover, $F$ is contained
in $K(x)^{alg}$, and therefore $F$ is a field. 
So,  $J_1$ is a maximal ideal, and  $J_1 \cap K[x] = (0)$. 
Indeed, the restriction of the projection $\pi: S  \rightarrow S /J_1$ to $K[x]$  
is injective since $K[ x] \subseteq K(x) \subseteq  S/J_1$. 
So, $\ker(\pi{\upharpoonright_{K[x]}}) = J _1\cap K[x] = (0)$, and hence $J_1$ is an
E-ideal of~$S$. 
Now we want to extend $J_1$ to a maximal E-ideal $J$ of  $K[ x]^E$. 
By repeating the proof of Lemma~\ref{lem:R-S} we have $S^E = K[ x]^E$. 
By Lemma \ref{ideali} there exists an exponential maximal ideal $J$ of $K[x]^E$ extending $J_1$ which is clearly not of the form $I_a$ for any $a\in K$.
%
\end{proof}



Using a similar construction and the ideas in the previous proposition, we
obtain the following more general result, which does not require that $K$ is
algebraically closed.

\begin{thm}\label{thm:max-zero}
 Let $K$ be an E-field and fix $n \ge 1$. 
There exists an E-ideal $J$  of $K[\overline x]^E$ such that:
\begin{enumerate}
\item $J$ is strongly maximal;
\item $J \cap R_n = (0)$.
\end{enumerate}
\end{thm}

\begin{proof} We construct an E-ideal $J_{n+1}$ of $R_{n+1}$ which
satisfies (1) and (2) of the statement.
The conclusion then follows from Lemma~\ref{ideali}.


\smallskip

We have $R_{n+1} = R_{n}[t^{A}]$ as $R_{n}$-algebra, for some suitable divisible
subgroup $A$ of~$R_{n}$. 
Fix
 a $\Q$-basis $\set{a_{i}: i \in I}$ of~$A$.
Let $\set{p_{i}: i \in I}$ be the set of nonzero elements of $R_{n}$ (we are using the same
index set~$I$ since $\card I = \card{R_n}$).
Let $F$ be the fraction field of $R_{n}$, and $\bar F$ be its algebraic
closure.
Since the multiplicative group $\bar F^{*}$ of $\bar F$ is divisible, there exists some (non
unique!) group homomorphism $\alpha: \tuple{A, +} \to \tuple{\bar F^{*}, \cdot}$, such that, for every $i \in I$,
$\alpha(a_{i}) = p_{i}$.
By the universal property of the group ring $R_{n}[t^{A}]$, there exists a
unique homomorphism of $R_{n}$-algebras $\beta: R_{n}[t^{A}] \to \bar F$  such that,
for every $a \in A$, $\beta(t^{a}) = \alpha(a)$.

Let $J_{n+1}$ be the kernel of $\beta$, and $L := R_{n+1}/J_{n+1}$.
Notice that $L$ is a subring of the field $\bar F,$ therefore $L$ is an integral
domain, and thus $J_{n+1}$ is a prime ideal.
Let $r \in J_{n+1} \cap R_n,$ then
\[
0 = \beta(r) = \beta(r \cdot 1) = r \cdot \beta(1).
\]
Since $\bar F$ is a field and $1 = \beta(1)$ we have  $r = 0$.
Thus, $J_{n+1} \cap R_{n} = (0)$.

Finally, $L \simeq \beta(R_{n+1})$, so $L$ is an algebraic extension of the field of fractions of
$R_{n}$, and therefore $L$ is a field, proving that $J_{n+1}$ is a maximal ideal.
\end{proof}

\section{Exponential kernel}\label{sec:Ekernel}
Let $R = (R; V,E)$ be a partial E-ring and
$I < R$ be  an E-ideal.\\
We define the E-kernel of $I$ as
\[
Z(I) := \set{a \in V: E(a) - 1 \in I},
\] 
and the exponential kernel of $R$ as
\[
\ker(E_{R}) := \set{a \in V: E(a) = 1}.
\]
Notice that $I \cap V \subseteq Z(I) \subseteq V$ and $Z(I)/I$ is the
exponential kernel of  the quotient $R/I$ (we use the notation $B/C$ for $B, C$
subgroups of the Abelian group $A$ to denote $(B+C)/C$): that is,
\[
Z(I)  =
\set{a \in V: E( a + I) = 1 \text{ inside } R/I}.
\]

\medskip
  
In Sections~\ref{sec:smideal} and~\ref {sec:smideals-counterexample}
we proved that the E-ideals of $K[\bar x]^{E}$ of the form $I_{\av}$ are maximal, but not
all maximal E-ideals are of such form.
Moreover, if $I = I_{\av}$, then
the exponential kernel does not increase, i.e.
$Z(I_{\av})/I_{\av} = \ker(E_{K})$. 
We ask
if the converse is true, i.e.\ we ask if (at least, when $K$ is algebraically closed),
given a maximal 
E-ideal of $K[x]^E$ such that $Z(I)/I = \ker(E_{K})$, we have that
$I = I_{a}$ for some~$a \in K$. 
Corollary~\ref{cor:E-kernel} gives counterexamples for \emph{prime} E-ideals:
what happens for maximal E-ideals is still open.

\begin{lemma}
Let $\bf R'$ be the free 1-extension of $R$ and $I' < \bf R'$ be the E-ideal of $\bf R'$
generated by~$I$. 
Then,
\begin{equation}\label{eq:Z-1}
Z(I') = Z(I) + I.
\end{equation}
Moreover, let $J$ be the E-ideal of $R^{E}$ generated by $I$.
Then, 
\begin{equation}\label{eq:Z-2}
Z(J) = Z(I) + J.
\end{equation}
\end{lemma}
\begin{proof}
Notice that $Z(I) < V$ and $Z(I') < R$.
Notice also that \eqref{eq:Z-2} follows by iteration from \eqref{eq:Z-1}.

Let $P < I$ and $Q < R$ be $\Q$-vector spaces such that $R = V \oplus P \oplus Q$.

Let $\sigma_{P}: R[t^{P}] \to R$ be the canonical augmentation map.
Notice that ${\bf R'} = R[t^{P}][t^{Q}]$.
Moreover, 
\[
I' = {\bf R'} \cdot \sigma_{P}^{-1}(I). 
\]
Let $a \in Z(I')$.
We want to show that $a \in Z(I) + I$.
We can write uniquely $a = v + p  + q$, with $v \in V, p \in P, q \in Q$.
\begin{claim*}
$q = 0$.
\end{claim*}
In fact, an element $b \in \bf R'$ can be written uniquely as
\[
b = \sum_{x \in Q} b_{x} E(x)
\]
where each $b_{x} \in R[t^{P}]$.
Such $b$ is in $I'$ iff each $b_{x} \in \sigma_{P}^{-1}(I)$.
If $q \ne 0$ and $b = E(a) - 1 \in I',$ then $b_{0} = -1 \notin \sigma_{P}^{-1}(I)$.

Thus, we have $a = v + p$.
Moreover, since $a \in Z(I')$, we have 
\[
I \ni \sigma_{P}(E(a)-1) = \sigma_{p}(E(v)) \sigma_{P}(E(p)) - 1
= E(v) - 1.
\] 
showing that $v \in Z(I)$.
Therefore, 
\[
a \in Z(I) + P \subseteq Z(I) + I.
\]

Conversely, assume that $a \in Z(I) + I$.
Write $a = v + b$, with $v \in Z(I)$ and $b \in I$.
We have
\[
E(a) = E(v) E(b) \in (1 + I) \cdot (1 + I) \subseteq  1 + I. \qedhere
\]
\end{proof}


\begin{corollary}\label{cor:E-kernel}
If $K$ is uncountable, 
there exists an E-ideal $I < K[x]^E$ such that:
\begin{enumerate}
\item $I$ is prime;
\item $Z(I) = \ker(E_{K}) + I$;
\item $x - a \notin I$ for every $a \in K$.
\end{enumerate}
\end{corollary}
It is easy to see that each of the above property is equivalent to the corresponding property:
\begin{enumerate}[(1')]
\item $L := K[x]^{E}/I$ is an integral domain (containing $K$);
\item $\ker(E_{L}) = \ker(E_{K})$;
\item $x + I \notin K$ (under the identification between $K$ and its image in~$L$).
\end{enumerate}
\begin{proof}
Let $J < {\bf R'}$ be an E-ideal satisfying:
\begin{enumerate}[(i)]
\item $J$ is a maximal ideal;
\item $Z(J)/J = \ker(E_{K})$;
\item $x - a \notin J$ for every $a \in K$.
\end{enumerate}
Let $I < K[x]^{E}$ be the E-ideal generated by $J$.
Then, by the above lemma, $I$ satisfies (1)--(3).

To produce $J$ as above, proceed as in the proof of Theorem~\ref{thm:max-nontrivial}:
fix a $\Q$-linear basis $B$ of $K$.
Let
\[
\set{y_{c} : c \in K} \subset {\bf R'}
\]
be an enumeration of $E(Bx) \subset {\bf R'}$.
Define
\[
J := \Pa{y_c - (x- c): c \in K}. 
\]
(the ideal generated by all the $y_{c} - (c-x)$ inside ${\bf R'}$).
Then, $J$ is an E-ideal satisfying the above conditions (i)--(iii).
\end{proof}

\section{Exponential radical ideals and characterization}

In this section we introduce the notion of E-radical ideal for any E-ring $R$ as follows.

\begin{definition} Let $J$ be an E-ideal of an E-ring $R$. 
We define the E-radical ideal of $J$ as $\Erad(J) :=  \bigcap_{P \supseteq J} P$, where $P$
varies among prime E-ideals. 
\end{definition}

We aim at connecting the radical of an E-ideal of a polynomial ring to the variety associated to an E-ideal as in the case of classical polynomial rings. 

Let $K$ be an E-field and $J$  an E-ideal of the E-polynomial ring $K[\x]^E$.
Let $F$ be an E-domain extending $K$ (often, $F$~will be an E-field).
 
We  define  $\mathcal{I} (V(J))$ as follows:\\
$V_F(J) := \{ \overline a \in F^n : f(\overline a) = 0 \mbox{ for all } f(\x) \in J
\}$,\\
$V(J) := \bigcup V_F(J)$ as $F$ varies among all E-fields containing $K$, and\\
$ \mathcal{I} (V(J)) = \{ p(\x) \in K[\x]^E : p(\overline a) = 0 \mbox{ for all } \overline
a \in V(J) \}$. 

\begin{remark}
Let $(R,E)$ be an E-domain and $F$ be its fraction field.
Then, there exists at least one way to extend the exponential function to all of $F$.
\end{remark}
\begin{proof}
Let $A\subset F$ be a complement of $R$ as $\Q$-linear spaces.
For every $r \in R$ and $a \in A$, define $E'(a + r) := E(r)$.
Then, $E'$ is an exponential function on $F$ extending~$E$.
\end{proof}

\begin{corollary}
$\mathcal{I} (V(J)) = 
\{ p(\x) \in K[\x]^E : p(\overline a) = 0 \mbox{ for all } \overline a \in V_{F}(J)$ 
as $F$ varies among all E-domains containing $K\}$. 
\end{corollary}

\begin{lemma}
Let  $J$ be an E-ideal of the  E-polynomial ring $K[\x]^E.$ Then $\Erad(J) = \mathcal{I} (V(J)) $.
 \end{lemma}

\begin{remark}\label{idealeradicale}
The E-ideal $I = (xy)^E$ is not an E-radical ideal, unlike in the classical
case. 
Indeed $I$ is not the intersection of prime E-ideals, since $(xy)^E \neq (x)^E
\cap (y)^E,$ because $x(E(y) - 1) \in (x)^E \cap (y)^E \setminus (xy)^E$. 
\end{remark}

The situation is simpler for a certain subring of $\C[x]^E$ where two polynomials with the same roots generated the same radical ideal (see \cite{FT-poly}).

\begin{remark}
If $J$ is not contained in any prime E-ideals then $\Erad(J) = K[\x]^E.$  Let  $I = (xy, E(x) + 1, E(y) + 1)$ be an ideal of $S = K[x, y, e^{\Q x}, e^{\Q y}]$ where $S$ is a subring of $K[\x]^E.$ In particular,  $I$  is an E-ideal of $S$. Then 
$I^E$ is  an E-ideal of $K[\x]^E$ which is not contained in any
prime E-ideal. Indeed, if $I \subseteq P$ where $P$ is prime E-ideal, then $xy \in I \subseteq P.$ 
Since $P$ is  prime,  \wloG   $x \in P$, and so $E(x) - 1 \in P$.
But $E(x) + 1 \in I \subseteq P$, and thus $2 \in P$,  contradiction. 
\end{remark}

We consider an  E-ring $R$ and let  $J$ be an E-ideal of~$R$.  

We study prime E-ideals and  E-radical ideals, i.e. E-ideals which are equal to
their E-radical.
We characterize $\Erad(J)$ using the following theory.\\

We consider the following first-order language 
$\mathcal{L} = \{+, -, \cdot, e^{x}, 0, 1\} \cup \{ S \}$ where $S$ is a unary predicate.\\

We recall the following definition of Horn clause which we will use  in analyzing classes of structures as in the context of universal algebra.


\begin{definition}
Given $p_1, \ldots, p_k, q$ $\mathcal{L} $-terms,  i.e.\ exponential polynomials with integers coefficients, we define the associated strict Horn clauses the formulas of the form
$$S(p_1) \wedge \ldots \wedge S(p_k) \rightarrow S(q),$$ 
in another words if 
\[ p_1 \in S \wedge \ldots \wedge p_k \in S \mbox{ then } q \in S.\]

\end{definition}


\begin{examples}
$S(x^2)  \rightarrow S(x)$, $S(x) \wedge S(y)  \rightarrow S(x + y) ,$ $S(x)  \rightarrow S(e^x - 1)$ and $S(0)$ are $\mathcal{L} $-strict Horn clause.
\end{examples}

In order to lighten the notation, we write $p_1 \wedge \ldots \wedge p_k \rightarrow q$ in place of
 $S(p_1) \wedge \ldots \wedge S(p_k) \rightarrow S(q)$.

\smallskip

\begin{notation} In the following we use $(R, J) \models \alpha,$ where $R$ is an E-ring, $J$ is an ideal, or an E-ideal or a prime E-ideal, to say that $(R, J)$ satisfies a Horn strict clause $\alpha,$ i.e. for any $a \in R$ if $p_1(a) \wedge \ldots \wedge p_k(a) \in J$ then $q(a) \in J.$ 
\end{notation}

We consider the following theories:\\
\[\begin{aligned}
T_{0} &:= \{ \mbox{ Horn clause } \alpha : \mbox{ for any E-ring R and ideal J, } (R, J) \models \alpha \}\\
T_{1} &:= \{ \mbox{ Horn clause } \alpha : \mbox{ for any E-ring R and E-ideal J, } (R, J) \models \alpha \}\\
T_p &:=  \{ \mbox{ Horn clause } \alpha : \mbox{ for any E-ring R and prime E-ideal
  J, } (R, J) \models \alpha \}.
\end{aligned}
\]

Clearly, $T_{1}$ is generated by the following Horn clauses:
\[ 0, \quad x \wedge y  \rightarrow x-y, \quad x  \rightarrow xy, \quad x \rightarrow e^{x}-1.\]

Our aim is to give an explicit description of $T_{p}$ and relate it to the E-radical.
In $T_p$ we have, besides the clauses in $T_{1}$, also others; for instance,
the following clauses are in $T_{1}$:
$x^n  \rightarrow x$,  $(xy \wedge e^{x} + 1)  \rightarrow y$.

Let $T$ be  a set of Horn clauses.
Let $\model{T}_{R}$ be the  following family   of subsets of $R$:
\[
\model{T}_{R}  := \set{J \subseteq R: (R, J) \models T}. 
\]
In what follows we will drop the subscript $R$.
\begin{remark}\label{unione e intersezione} 
$\model{T}$ is closed under arbitrary intersections and 
under increasing unions.
\end{remark}
Thus, we can consider the ``radical'' operator associated to $T$.
We will let $X$ vary among subsets of~$R$.

We define $\Trad X$ as the intersection of the sets in $\model T$ 
containing~$X$.
We have that $\Trad X$ is the smallest subset of $R$ containing $X$ and such
that $(R, \Trad X ) \models T$.
In particular, $(R, X) \models T$ iff $X = \Trad X$.

We can build $\Trad J$ in a  ``constructive'' way,  \ie we have a description
of all elements of $\Trad J$.

\begin{definition}
Given a family $\Ffam$ of 
$\mathcal{L}$-structures, we denote its theory by
\[
Th(\Ffam) := \set{\alpha: \alpha \mbox{ is a Horn clause and for all }  M \in \Ffam, \ M \models \alpha }.
\]
The deductive closure of $T$ is $\Tbar := Th(\model T)$.
We say that  $T$ clauses is deductively closed if $T = Th(\Ffam)$ for some
family $\Ffam$, or equivalently if $T = \Tbar$.
An axiomatization of $T$ is a set of Horn clauses $S$ such that
$\overline S = \Tbar$.
\end{definition}

\begin{example}

\begin{enumerate}

\item $\Tunorad X = (X)^{E}$, the E-ideal generated by~$X$.

\item Let $T'_{0} := \set{0, x \wedge y \rightarrow x - y, x  \rightarrow xy}$; notice that $T_0$ it is the deductive closure of $T'_0.$
Then, $\model{T'_{0}} = \model{T_{0}}$ is the family of pairs $(R,J)$ with $J$ ideal of $R$.
We have that 
$\Tzprimerad$
is the ideal generated by~$X$.

\item Let $T_{R} := \set{0, x \wedge y \rightarrow x - y, x  \rightarrow xy, x^{2} \rightarrow x}$.
Then, $\model{T_{R}}$ is the family of pairs $(R,J)$ with $J$ radical ideal of $R$.
We have that $\radicale{T_{R}}(X)$ is the radical of the ideal generated by~$X$.

\end{enumerate}
\end{example}

\begin{remark}
\[
\Trad X  = \set{q(\cv): 
p_{1}(\x) \wedge \dots \wedge p_{k}(\x) \rightarrow q(\x) \in \Tbar,
 \cv \in R^{< \omega}, p_i(\cv) \in X, i = 1, \dotsc, k}. 
\]
\end{remark}
\begin{proof}
Let us denote 
$Y := \set{q(\cv): 
p_{1}(\x) \wedge \dots \wedge p_{k}(\x) \rightarrow q(\x) \in \Tbar,
 \cv \in R^{< \omega}, p_i(\cv) \in X, i = 1, \dotsc, k}$.

It is clear that $X \subseteq Y \subseteq \Trad X$.
We want to show that $\Trad X \subseteq Y$.
It suffices to show that $Y \in \model T$.
Let $p_{1}(\x) \wedge \dots \wedge p_{k}(\x) \rightarrow q(\x) \in T$ and $\cv \in R^{< \omega}$ such
that
$p_{i}(\cv) \in Y$, $i = 1, \dotsc, k$. 
It suffices to show the following:
\begin{claim}
$q(\cv) \in Y$.
\end{claim}
For simplicity of notation, we assume that  $k = 1$ and $p := p_{1}$.
Since $p(\cv) \in Y$, by definition of~$Y$, there exist $r_{1}(\x, \x') \wedge
\dots \wedge r_{\ell}(\x, \x') \rightarrow p(\x) \in \bar T$ and $\cv' \in R^{< \omega}$
such that $r_{j}(\cv,\cv') \in X$, $i = 1, \dotsc, \ell$.
Notice that every $M \in \model T$ satisfies
\[
\beta := \bigwedge_{\ell} r_{\ell}(\x, \x') \rightarrow q(\x).
\]
Therefore, $\beta \in  \bar T$, and therefore, by definition,
$q(\cv) \in Y$.
\end{proof}

\begin{corollary}\label{cor:Trad-union}
For every $b \in R$,
\begin{multline*}
\Trad{Xb} = \set{q(\cv): 
p_{1}(\x) \wedge \dots \wedge p_{k}(\x) \rightarrow q(\x) \in \Tbar, \\ 
 \cv \in R^{< \omega}, \ p_i(\cv) \in X \vee p_{i}(\cv) = b, \ i = 1, \dotsc, k}. 
\end{multline*}
\end{corollary}

\smallskip

Now we want to give more explicit description of $ \Tzrad  X $. In order to do it we define a new operator $\sqrt[E] X$ on subsets of $R$.\\


For every $n \in N$, we define the operator $\sqrt[n]{\vphantom{x} }$ on subsets of $R$
inductively in the following way:\\
\[\begin{aligned}
\sqrt[0] X &= (X)^E;\\
\sqrt[1] X &=  \sqrt[0] { \{ a \in R:  
\exists b_1, b_2 \in R : b_1 \cdot b_2  \in
  \sqrt[0] X  \wedge a  \in \sqrt[0]{X b_1} \cap \sqrt[0]{X b_2} \} };\\
\sqrt[n+1] X &= \sqrt[n] { \{ a :  
\exists
b_1, b_2  : b_1 \cdot b_2 \in  \sqrt[n] X  \wedge a  \in \sqrt[n] {X b_1} \cap \sqrt[n]{X b_2} \} }.
\end{aligned}
\]

Let $\sqrt[E] X := \bigcup_{n \in \N} \sqrt[n] X.$

\begin{thm}\label{thm:Erad-up}
$\Erad X = \Tzrad X = \sqrt[E] X$.
\end{thm}

\proof
In the proof we need the following lemmas.
\begin{lemma}\label{lem:Erad-prod}
Let $J \subseteq R$.
Let $b_{1} \cdot b_{2} \in J$.
Assume that $a \in \Tzrad{J b_{1}} \cap \Tzrad{J b_{2}}$.
Then, $a \in \Tzrad J$.
\end{lemma}
\begin{proof}
If $J$ were a prime E-ideal, the result would be clear.
In general, by Corollary~\ref{cor:Trad-union}, there exist Horn clauses
$p_{i,1}(\x) \wedge \dots \wedge p_{i,k}(\x) \wedge z \rightarrow q_{i}(z,\x) \in T_{p}$, $i = 1,2$,
and
$\cv \in R^{< \omega}$, such that
\[
a = q_{i}(b_{i},\cv), 
\quad
p_{i,j}(\cv) \in J,
\quad i = 1,2, 
\quad j = 1, \dotsc, k.
\]
Moreover, the Horn clause
\[
\alpha(w, z_{1}, z_{2}, \x)  := z_{1} \cdot z_{2} \wedge w - q_{1}(z_{1}, \x) \wedge w - q_{2}(z_{2},\x)
\wedge \bigwedge_{i,k} p_{i,k}(\x) 
\rightarrow w
\]
is in $T_{p}$ (since it is satisfied by any prime E-ideal).
Thus, $(R, \Tzrad J) \models \alpha$ and the conclusion follows by considering
$\alpha(a, b_{1}, b_{2}, \cv)$.
\end{proof}

Therefore, $\Erad X \supseteq \Tzrad X \supseteq \sqrt[E] X$.
Thus, it suffices to show that $\Erad X \subseteq \sqrt[E] X$, or, equivalently, that
for every $a \in R \setminus \sqrt[E]X$, we have $a \notin \Erad X$.

 



We need a further result.

\begin{lemma}\label{lem:Trad-max}
Let  $a \in R$ and $P$ be an E-ideal of $R$ maximal among the E-ideals $J $of $R$
not  containing~$a$ and such that $\sqrt[E] J = J$.
Then, $P$ is prime.
\end{lemma}

\begin{proof}
Suppose $b_1 \cdot b_2 \in P$ and by contradiction we assume that $b_1, b_2 \not \in P$.
We define $Q_1 = \sqrt[E]{P  b_1}$ and $Q_2 = \sqrt[E]{P  b_2}$. 
We have that  $Q_{1},Q_{2} \in \model{T_{0}}$.
Notice that $Q_1, Q_2 \supset P$, and therefore the maximality of $P$ implies that $a
\in Q_{1}$ and $a \in Q_{2}$.
Moreover, $b_1 \cdot b_2 \in P$; since $\sqrt[E] P = P$, we have
  that $a \in  \Tzrad {P } = P$,  contradiction.   
\end{proof}

We can now conclude the proof of Theorem.~\ref{thm:Erad-up}.
Let  $a \notin \sqrt[E]X$.
By Lemma~\ref{lem:Trad-max}, there exists a prime E-ideal $P$
containing $X$ such that $a \notin P$.
Thus, $a \notin \Erad X$, and we are done.
\qed

\begin{thm}\label{thm:Erad}
$J$ is an E-radical iff, for every $a, b_{1}, b_{2} \in R$
\[
b_{1} \cdot b_{2} \in J \wedge a \in \sqrt[E]{J b_{1}} \cap \sqrt[E]{J b_{1}} 
\rightarrow a \in J.
\]
\end{thm}
From the above Theorem we can extract a recursive axiomatization of $T_{p}$, 
using
Corollary~\ref{cor:Trad-union} to characterize $\sqrt[E]{J b_{i}}$.

\bigskip

We can interpret the discussion above in terms of quasi-varieties.
\begin{definition}
An E-ring is E-reduced if $\Erad(0) = (0)$.
\end{definition}
By Theorem~\ref{thm:Erad-up}, the class E-red
of E-reduced E-rings is a quasi-variety:
\ie, it can be axiomatized via Horn formulae in the language
\[
\{ + , \cdot, - , E, 0, 1\},
\]
where $+, \cdot$ are binary operations, $-, E$ are unary operations, and $0,1$ are constants
(we replace the condition $t \in (0)$ with the condition $t = 0$).
By \cite{Burris-Sankappanavar}*{Thm.~2.25}, E-red is closed under isomorphisms, taking substructures, reduced product.\\
It is not difficult to see directly that E-red is closed under isomorphisms,  substructures, direct
products,  and ultraproducts. Thus, by \cite{Burris-Sankappanavar}*{Thm.~2.25},
E-red is a quasi-variety.
From this it is easy to deduce that $\Erad = \Tzradop$.
With the proof we gave of Theorem~\ref{thm:Erad-up} we obtained the extra
information that $\Erad = \sqrt[E]{\vphantom{x}}$, and Theorem~\ref{thm:Erad},
with the recursive axiomatization of~$T_{p}$.

Moreover, Lemma~\ref{lem:Trad-max} might be of independent interest.


\bibliography{exponential}		

\begin{bibdiv}
\begin{biblist}

\bib{Arsla}{article}{
      author={Aslanyan, Vahagn},
      author={Kirby, Jonathan},
      author={Mantova, Vincenzo},
       title={A geometric approach to some systems of exponential equations},
        date={2023},
     journal={International Mathematics Research Notices},
      volume={2023},
      number={5},
       pages={4046\ndash 4081},
}

\bib{bayskirby}{article}{
      author={Bays, Martin},
      author={Kirby, Jonathan},
       title={Pseudo-exponential maps, variants, and quasiminimality},
        date={2018},
        ISSN={1937-0652},
     journal={Algebra Number Theory},
      volume={12},
      number={3},
       pages={493\ndash 549},
         url={https://doi.org/10.2140/ant.2018.12.493},
      review={\MR{3815305}},
}

\bib{bourbaki}{book}{
      author={Bourbaki, Nicolas},
       title={Algebra {II}},
   publisher={Springer},
        date={1981},
}

\bib{Burris-Sankappanavar}{book}{
      author={Burris, Stanley},
      author={Sankappanavar, H.~P.},
       title={A course in universal algebra},
    language={English},
      series={Grad. Texts Math.},
   publisher={Springer, Cham},
        date={1981},
      volume={78},
}

\bib{DFT1}{article}{
      author={D'Aquino, Paola},
      author={Fornasiero, Antongiulio},
      author={Terzo, Giuseppina},
       title={Generic solutions of equations with iterated exponentials},
        date={2018},
        ISSN={0002-9947},
     journal={Trans. Amer. Math. Soc.},
      volume={370},
      number={2},
       pages={1393\ndash 1407},
         url={https://doi.org/10.1090/tran/7206},
      review={\MR{3729505}},
}

\bib{DFT2}{article}{
      author={D'Aquino, Paola},
      author={Fornasiero, Antongiulio},
      author={Terzo, Giuseppina},
       title={A weak version of the strong exponential closure},
        date={2021},
        ISSN={0021-2172},
     journal={Israel J. Math.},
      volume={242},
      number={2},
       pages={697\ndash 705},
         url={https://doi.org/10.1007/s11856-021-2141-1},
      review={\MR{4282096}},
}

\bib{null}{article}{
      author={D'Aquino, Paola},
      author={Macintyre, Angus},
      author={Terzo, Giuseppina},
       title={Schanuel {N}ullstellensatz for {Z}ilber fields},
        date={2010},
        ISSN={0016-2736},
     journal={Fund. Math.},
      volume={207},
      number={2},
       pages={123\ndash 143},
         url={https://doi.org/10.4064/fm207-2-2},
      review={\MR{2586007}},
}

\bib{shapiro}{article}{
      author={D'Aquino, Paola},
      author={Macintyre, Angus},
      author={Terzo, Giuseppina},
       title={From {S}chanuel's conjecture to {S}hapiro's conjecture},
        date={2014},
        ISSN={0010-2571},
     journal={Comment. Math. Helv.},
      volume={89},
      number={3},
       pages={597\ndash 616},
         url={https://doi.org/10.4171/CMH/328},
      review={\MR{3260843}},
}

\bib{Everest_Porten}{article}{
      author={Everest, G.~R.},
      author={van~der Poorten, A.~J.},
       title={Factorisation in the ring of exponential polynomials},
        date={1997},
        ISSN={0002-9939},
     journal={Proc. Amer. Math. Soc.},
      volume={125},
      number={5},
       pages={1293\ndash 1298},
         url={https://doi.org/10.1090/S0002-9939-97-03919-1},
      review={\MR{1401740}},
}

\bib{FT-poly}{article}{
      author={Fornasiero, Antongiulio},
      author={Terzo, Giuseppina},
       title={A note on exponential polynomials},
        date={2024-07-04},
        ISSN={1827-3491},
     journal={Ricerche di Matematica},
      eprint={https://doi.org/10.1007/s11587-024-00871-8},
         url={https://doi.org/10.1007/s11587-024-00871-8},
}

\bib{Gall}{article}{
      author={Gallinaro, Francesco~Paolo},
       title={Exponential sums equations and tropical geometry},
        date={2023},
     journal={Selecta Mathematica, New Series},
      volume={29},
      number={4},
}

\bib{GP74}{article}{
      author={Gilmer, Robert},
      author={Parker, Tom},
       title={Divisibility properties in semigroup rings},
        date={1974},
        ISSN={0026-2285},
     journal={Michigan Math. J.},
      volume={21},
       pages={65\ndash 86},
         url={http://projecteuclid.org/euclid.mmj/1029001210},
      review={\MR{342635}},
}

\bib{HR}{article}{
      author={Henson, C.~Ward},
      author={Rubel, Lee~A.},
       title={Some applications of {N}evanlinna theory to mathematical logic:
  identities of exponential functions},
        date={1984},
        ISSN={0002-9947},
     journal={Trans. Amer. Math. Soc.},
      volume={282},
      number={1},
       pages={1\ndash 32},
         url={https://doi.org/10.2307/1999575},
      review={\MR{728700}},
}

\bib{HRS}{article}{
      author={Henson, C.~Ward},
      author={Rubel, Lee~A.},
      author={Singer, Michael~F.},
       title={Algebraic properties of the ring of general exponential
  polynomials},
        date={1989},
        ISSN={0278-1077},
     journal={Complex Variables Theory Appl.},
      volume={13},
      number={1-2},
       pages={1\ndash 20},
         url={https://doi.org/10.1080/17476938908814374},
      review={\MR{1029352}},
}

\bib{Katzberg}{article}{
      author={Katzberg, H.},
       title={Complex exponential terms with only finitely many zeros},
        date={1983},
      volume={49},
       pages={68–72},
}

\bib{Kolchin}{book}{
      author={Kolchin, E.~R.},
       title={Differential algebra and algebraic groups},
      series={Pure and Applied Mathematics, Vol. 54},
   publisher={Academic Press, New York-London},
        date={1973},
      review={\MR{568864}},
}

\bib{MAC_free}{article}{
      author={Macintyre, Angus},
       title={Schanuel's conjecture and free exponential rings},
        date={1991},
        ISSN={0168-0072},
     journal={Ann. Pure Appl. Logic},
      volume={51},
      number={3},
       pages={241\ndash 246},
         url={https://doi.org/10.1016/0168-0072(91)90017-G},
      review={\MR{1098783}},
}

\bib{macintyre}{incollection}{
      author={Macintyre, Angus},
       title={Exponential algebra},
        date={1996},
   booktitle={Logic and algebra ({P}ontignano, 1994)},
      series={Lecture Notes in Pure and Appl. Math.},
      volume={180},
   publisher={Dekker, New York},
       pages={191\ndash 210},
      review={\MR{1404940}},
}

\bib{macWilkie}{incollection}{
      author={Macintyre, Angus},
      author={Wilkie, A.~J.},
       title={On the decidability of the real exponential field},
        date={1996},
   booktitle={Kreiseliana},
   publisher={A K Peters, Wellesley, MA},
       pages={441\ndash 467},
      review={\MR{1435773}},
}

\bib{Manders}{article}{
      author={Manders, Kenneth},
       title={On algebraic geometry over rings with exponentiation},
        date={1987},
        ISSN={0044-3050},
     journal={Z. Math. Logik Grundlag. Math.},
      volume={33},
      number={4},
       pages={289\ndash 292},
         url={https://doi.org/10.1002/malq.19870330402},
      review={\MR{913290}},
}

\bib{mantova}{article}{
      author={Mantova, Vincenzo},
       title={Polynomial-exponential equations and {Z}ilber's conjecture},
        date={2016},
        ISSN={0024-6093},
     journal={Bull. Lond. Math. Soc.},
      volume={48},
      number={2},
       pages={309\ndash 320},
         url={https://doi.org/10.1112/blms/bdv096},
        note={With an appendix by Mantova and U. Zannier},
      review={\MR{3483068}},
}

\bib{marker}{article}{
      author={Marker, David},
       title={A remark on {Z}ilber's pseudoexponentiation},
        date={2006},
        ISSN={0022-4812},
     journal={J. Symbolic Logic},
      volume={71},
      number={3},
       pages={791\ndash 798},
         url={https://doi.org/10.2178/jsl/1154698577},
      review={\MR{2250821}},
}

\bib{Matsumura}{book}{
      author={Matsumura, Hideyuki},
       title={Commutative algebra},
   publisher={W. A. Benjamin, Inc., New York},
        date={1970},
      review={\MR{0266911}},
}

\bib{point}{unpublished}{
      author={Point, F.},
      author={Regnault, N.},
       title={Exponential ideals and a nullstellensatz},
        date={2020},
        note={Preprint},
}

\bib{Ritt}{article}{
      author={Ritt, J.~F.},
       title={A factorization theory for functions
  {$\sum_{i=1}^na_ie^{\alpha_ix}$}},
        date={1927},
        ISSN={0002-9947},
     journal={Trans. Amer. Math. Soc.},
      volume={29},
      number={3},
       pages={584\ndash 596},
         url={https://doi.org/10.2307/1989097},
      review={\MR{1501406}},
}

\bib{terzo}{article}{
      author={Terzo, Giuseppina},
       title={Some consequences of {S}chanuel's conjecture in exponential
  rings},
        date={2008},
        ISSN={0092-7872},
     journal={Comm. Algebra},
      volume={36},
      number={3},
       pages={1171\ndash 1189},
         url={https://doi.org/10.1080/00927870701410694},
      review={\MR{2394281}},
}

\bib{tougeron}{article}{
      author={Tougeron, Jean-Claude},
       title={Alg\`ebres analytiques topologiquement noeth\'{e}riennes.
  {T}h\'{e}orie de {K}hovanski\u{\i}},
        date={1991},
        ISSN={0373-0956},
     journal={Ann. Inst. Fourier (Grenoble)},
      volume={41},
      number={4},
       pages={823\ndash 840},
         url={http://www.numdam.org/item?id=AIF_1991__41_4_823_0},
      review={\MR{1150568}},
}

\bib{van}{article}{
      author={van~den Dries, Lou},
       title={Exponential rings, exponential polynomials and exponential
  functions},
        date={1984},
        ISSN={0030-8730},
     journal={Pacific J. Math.},
      volume={113},
      number={1},
       pages={51\ndash 66},
         url={http://projecteuclid.org/euclid.pjm/1102709376},
      review={\MR{745594}},
}

\end{biblist}
\end{bibdiv}

\end{document}